\documentclass[11pt]{article}
\usepackage{amsmath,amsfonts,epsfig}
\newcommand{\bb}{\mathbb}

\newcommand{\cx}{{\bb C}}

\newcommand{\integers}{{\bb Z}}

\newcommand{\reals}{{\bb R}}

\newcommand{\hthree}{{\bb H}^3}

\newcommand{\ddz}{\frac{\partial}{\partial z}}
\newcommand{\widemargins}{
\setlength{\textwidth}{6.0in}
\setlength{\oddsidemargin}{0.25in}
\setlength{\evensidemargin}{0.25in}
}

\newcommand{\qed}[1]{\nopagebreak[4]{\tiny \hfill
\fbox{\ref{#1}} \linebreak }\pagebreak[2]}
\newcommand{\del}{\partial}

\newcommand{\isom}{\cong}

\newcommand{\zbar}{{\overline{z}}}

\newcommand{\chat}{\widehat{\cx}}

\newcommand{\ad}{\operatorname{ad}}

\newcommand{\area}{\operatorname{area}}

\newcommand{\inj}{\operatorname{inj}}

\renewcommand{\Im}{\operatorname{Im}}

\renewcommand{\Re}{\operatorname{Re}}

\newcommand{\sech}{\operatorname{sech}}

\newcommand{\str}{\operatorname{str}}

\newcommand{\tr}{\operatorname{tr}}
\newcommand{\vol}{\operatorname{vol}}
\newcommand{\curl}{\operatorname{curl}}

\newtheorem{theorem}{Theorem}[section]
\newtheorem{prop}[theorem]{Proposition}
\newtheorem{lemma}[theorem]{Lemma}
\newtheorem{cor}[theorem]{Corollary}

\newcommand{\cC}{{\cal C}}
\newcommand{\cD}{{\cal D}}
\newcommand{\cE}{{\cal E}}
\newcommand{\cF}{{\cal F}}
\newcommand{\cG}{{\cal G}}
\newcommand{\cH}{{\cal H}}

\newcommand{\cL}{{\cal L}}
\newcommand{\cM}{{\cal M}}

\newcommand{\cR}{{\cal R}}
\newcommand{\cS}{{\cal S}}

\newcommand{\cU}{{\cal U}}

\widemargins
\newcommand{\lip}{{\operatorname{lip}}}
\newcommand{\bilip}{{\operatorname{bilip}}}
\newcommand{\width}{\operatorname{width}}
\begin{document}
\title{Hyperbolic cone-manifolds, short geodesics and Schwarzian derivatives}

\author{K. Bromberg\footnote{Supported by a grant from the NSF}}

\date{November 5, 2002}

\maketitle

\begin{abstract}
\noindent
Given a geometrically finite hyperbolic cone-manifold, with the cone-singularity sufficiently short, we construct a one parameter family of cone-manifolds decreasing the cone angle to zero. We also control the geometry of this one parameter family via the Schwarzian derivative of the projective boundary and the length of closed geodesics.
\end{abstract}

With his hyperbolic Dehn surgery theorem and later the orbifold
theorem, Thurston demonstrated the power of using hyperbolic
cone-manifolds to understand complete, non-singular hyperbolic
3-manifolds. Hodgson and Kerckhoff introduced analytic techniques to 
the study of cone-manifolds that they have used to prove deep results
about finite volume hyperbolic 3-manifolds. In this paper we use
Hodgson and Kerckhoff's techniques to study infinite volume hyperbolic
3-manifolds. The results we will develop have many applications:
Bers' density conjecture, the density of cusps on the boundary of
quasi-conformal deformations spaces and for constructing type
preserving sequences of Kleinian groups.

The simplest example of the problem we will study is the
following: Let $M$ be a hyperbolic 3-manifold and $c$ a simple closed
geodesic in $M$. Then the topological manifold $M \backslash c$, also has a
complete hyperbolic metric which we call $\hat{M}$. How does the
geometry of $M$ compare to that of $\hat{M}$?  Before attempting to
answer such a question we need to note that if $M$ has infinite volume
the hyperbolic structure will not be unique. If we do not make
further restrictions on the choice of $\hat{M}$ then there is no
reason to expect that $M$ and $\hat{M}$ will be geometrically
close. If $M$ is convex co-compact there is a natural choice to make for $\hat{M}$. Namely $M$ is compactified by a conformal structure $X$. We then choose $\hat{M}$ to be the unique geometrically finite hyperbolic structure on $M \backslash c$ with conformal boundary $X$. 

We can now return to our question: How do the geometry of $M$ and $\hat{M}$ compare? We will quantify this question in two ways. We will measure the length of geodesics in $M$ and $\hat{M}$ and we will measure the geometry of the ends of $M$ and $\hat{M}$ by bounding the distance between the projective structures on their boundaries. What we will see is that the change in geometry is bounded by the length of the geodesic $c$ in the original manifold $M$.

Results of this type were first obtained by McMullen \cite{McMullen:cusps}, in the case of a quasifuchsian manifold, where the geodesic $c$ is also short on a component of the conformal boundary. This work has been extended to arbitrary geometrically finite manifolds by Canary, Culler, Hersonsky and Shalen \cite{Canary:Culler:Hersonsky:Shalen}. Their techniques are entirely different then ours and one goal of this paper is to give new proofs of their estimates. These estimates are a key step in proving the density of cusps in the boundary of quasiconformal deformation spaces of hyperbolic 3-manifolds.

Our original motivation was another application, the density conjecture. In \cite{Bromberg:projective} we construct a family of quasifuchsian cone-manifolds, with cone angle $4\pi$, approximating any singly degenerate hyperbolic manifold with arbitrarily short geodesics. The estimates developed here imply that the quasifuchsian cone-manifolds are geometrically close to smooth quasifuchsian manifolds which will approximate the degenerate manifold. This is a special case of the Bers' density conjecture.

In joint work with Brock \cite{Brock:Bromberg:density} we extend the results of this paper. In particular we show that there is a diffeomorphism from $M\backslash c$ to $\hat{M}$ that is bi-Lipschitz outside of a tubular neighborhood of $c$ where the bi-Lipschitz constant is bounded by a constant depending only on the length of the $c$. Using this result we are able to prove the density conjecture for all freely indecomposable Kleinian groups without parabolics. 

Another question is the following: Assume $\Gamma$ is a Kleinian group and $\Gamma_i$ is a sequence of geometrically finite Kleinian groups such that $\Gamma_i \rightarrow \Gamma$, algebraically. Does there exist a type preserving sequence $\Gamma'_i$ of geometrically finite groups also converging to $\Gamma$? Here type preserving means that if elements $\gamma_i$ converge to $\gamma$ then $\gamma$ is parabolic if and only if the $\gamma_i$ are parabolic. In joint work with Brock, Evans and Souto \cite{Brock:Bromberg:Evans:Souto} we show that the answer to the question is yes and the type preserving sequence can be constructed by pinching the short geodesics in the $\Gamma_i$ to cusps. The estimates developed here and extended in \cite{Brock:Bromberg:density} can then be used to show that the new sequence has the same limit. This question is important because in many cases the work of Anderson-Canary \cite{Anderson:Canary:cores, Anderson:Canary:corestwo} implies that type preserving sequences are strong. Work of Evans \cite{Evans:tameness}, expanding on work of Canary-Minsky \cite{Canary:Minsky:tameness} then implies that the limit is tame.

The starting point for this paper is the local parameterization of hyperbolic cone-manifolds developed by Hodgson and Kerckhoff for finite volume manifolds in \cite{Hodgson:Kerckhoff:cone} and \cite{Hodgson:Kerckhoff:tube} and extended to geometrically finite cone-manifolds in \cite{Bromberg:rigidity}. These local results tell us that we can make a small decrease in the cone angle. To decrease the cone angle to zero we need to ensure that there are no degenerations. There are three possible types of degenerations that need to be avoided. First, we need to make sure that the cone-singularity does not develop a point of self intersection. Second we must show that there is a lower bound for the injectivity radius of the cone-manifold. Finally the geometry of the geometrically finite ends must be controlled. The first two problems are taken care of by work of Hodgson and Kerckhoff. It is the last problem that is the main work of this paper.

We now outline the remainder of the paper.

In \S1 we make enough definitions to state our main theorems.

In \S2 we review background material on projective structures, hyperbolic cone-manifolds and deformations thereof.

In \S3 we describe some of Hodgson and Kerckhoff's results on tubes in cone-manifolds. In particular, they derive estimates on the radius of embedded tubes about the cone-singularity. We also use these ideas to show that embedded hyperbolic half spaces are disjoint from these tubular neighborhoods of the cone-singularity.

In \S4 we use the analytic deformation theory of cone-manifolds to control the length of geodesics as the cone angle decreases. Again the key estimates are work of Hodgson and Kerckhoff.

\S5 is the heart of the paper. In this section we show that the $L^2$-norm of the cone-manifold deformation bounds the change in projective structure. We first do this for a hyperbolic half space and then use the work from \S3 to find a large embedded half space in each geometrically finite end which allows us to globally bound the deformation of the entire projective structure.

The final step is to understand the geometric limit of  a sequence of cone-manifolds. Our approach is essentially the same as in \cite{Hodgson:Kerckhoff:dehn} although a bit of extra care is need to take care of the geometrically finite ends.

Although all the results of this paper hold for manifolds with rank two cusps they are, as is often the case, an annoyance and distraction from the main line of argument. For this reason we defer the discussion of rank two cusps until \S7 where we show how the $L^2$-norm bounds control the shape of the cusp. We then outline how this result can be used to finish the proof of the main theorems for manifolds with rank two cusps.

In \S8 we derive the estimates of McMullen and Canary et. al. mentioned above.

In the appendix we recount mean value theorems for harmonic vector and strain fields that Hodgson and Kerckhoff proved in an early version of \cite{Hodgson:Kerckhoff:dehn}.

{\bf Acknowledgments.} As should be obvious to the reader of the above outline, none of the results in this paper would be possible without the important work of Hodgson and Kerckhoff on hyperbolic cone-manifolds. Their analysis of harmonic forms in a neighborhood of the cone-singularity underly all the estimates derived in this paper.

I'd also like to thank Jeff Brock who originally suggested to me that cone-manifolds could be used to measure the effect of drilling a short geodesics in hyperbolic  3-manifolds.

Much of the work in this paper was done while I was a post-doc at the University of Michigan. I would like to thank the department for its hospitality.
\section{Geometrically finite hyperbolic cone-manifolds}
Before we state the main theorems we need to define the main object of study, geometrically finite hyperbolic cone-manifolds. Let $N$ be a compact, differentiable, 3-manifold with boundary,
$\cC$ a collection of simple closed curves in the interior of
$N$ and let $M$ be the interior of $N \backslash \cC$. 
A {\em hyperbolic cone-metric} is a complete metric $g$ on the
interior of $N$ that restricts to a
Riemannian metric with constant sectional curvature equal to $-1$ on $M$;
i.e. a hyperbolic metric. At all components $c$ of $\cC$ the metric
will be singular; in cylindrical 
coordinates $(r, \theta, z)$ the metric will locally have the form
$$dr^2 + \sinh^2 r d\theta^2 + \cosh^2 r dz^2$$
where $\theta$ is measured modulo the cone-angle $\alpha$.
In this coordinate system the singular locus will be
identified with the $z$-axis. Note that the cone angle will be
constant along each component of the cone singularity.

A {\em complex projective structure} on a surface is an atlas of
charts to $\chat$ with transition maps M\"{o}bius transformations. A hyperbolic cone-metric $g$ is {\em geometrically finite without
  rank one cusps} if it extends to a projective structure on the
non-toroidal components $\del_0 N$ of $\del N$. To state this more
precisely we  recall that hyperbolic 3-space $\hthree$ is
naturally compactified by the  Riemann sphere $\chat$. Then $g$ is
geometrically finite if each $p \in \del_0 N$ has a neighborhood $V$
in $N$ and a map $\phi:V \longrightarrow \hthree \cup \chat$ such
that $\phi$ restricted $V \cap \operatorname{int} N$ is an isometry
into $\hthree$ and $\phi$ restricted to $V \cap \del_0 N$ defines an
atlas for a projective structure on $\del_0 N$. As we will not discuss
rank one cusps in this paper (except briefly in \S8) we simply refer to such metrics as
geometrically finite.

A projective structure on a surface $S$ also determines a conformal structure on $S$. Moreover, for a fixed conformal structure there will be many projective structures. We will often need to distinguish between the {\em projective boundary} and {\em conformal boundary} of a geometrically finite hyperbolic manifold.

If $\del N$ contains a torus $T$ the behavior near infinity is different. A neighborhood of $T$ in $M$ will be foliated by Euclidean tori of a fixed conformal class with area decreasing exponentially as the tori exit the end. Such a neighborhood is a {\em rank two cusp}. More explicitly every rank two cusp is the quotient of a subspace of $\hthree$ (in the upper half space model) of the form $\{(z,t) | t \geq t_0\}$
by parabolic isometries $z \mapsto z + 1$ and $z \mapsto z + \tau$ where $\Im \tau > 0$. Note that $\tau$ is the Teichm\"uller parameter of the tori that foliate the cusp.

As we mentioned in the introduction we will postpone discussion of rank two cusps whenever possible. However they cannot be completely avoided because as the cone angle limits to zero a tubular neighborhood of the cone singularity will limit geometrically to a cusp. For this reason it is natural to think of a rank two cusp as a cone singularity with cone angle zero. This is of particular importance for the local parameterization for cone-manifolds that we are about to state.

Let $\cG\cF(N,\cC)$ be the space of geometrically finite hyperbolic
cone-metrics on the interior of $N$ with cone singularity $\cC$ and
assume $\cG\cF(N, \cC)$ has the compact-$C^\infty$ topology. Then
$GF(N,\cC)$ will be $\cG\cF(N, \cC)$ modulo diffeomorphisms which are
isotopic to the identity and that fix each component of $\cC$. An equivalence class of metrics $[g] \in GF(N, \cC)$ assigns to each component of $\cC$ a cone angle and to $\del_0 N$ a conformal structure in the Teichm\"uller space $T(\del_0 N)$. If $\cC$ has $n$ components there is a map
$$\Psi: GF(N, \cC) \longrightarrow [0,\infty)^n \times T(\del_0 N)$$
with $\Psi([g]) = (\alpha_1, \dots, \alpha_n, X)$ where $\alpha_i$ is
the cone-angle of the $i$th component of $\cC$ and $X$ is the
conformal boundary. Note that as we discussed in the previous paragraph we have allowed the possibility for the cone angle to be zero in which case the cone singularity becomes a rank two cusp. We then have the following local parameterization:

\begin{theorem}[\cite{Hodgson:Kerckhoff:cone}, \cite{Hodgson:Kerckhoff:tube}, \cite{Bromberg:rigidity}]
\label{param}
Let $[g] \in GF(N, \cC)$ be a geometrically finite hyperbolic
cone-metric. Suppose the tube radius of the cone singularity is $\geq
\sinh^{-1}\left(\frac{1}{\sqrt{2}}\right)$. Then $\Psi$ is a local
homeomorphism at $[g]$.
\end{theorem}

{\bf Remark.} Hodgson and Kerckhoff first proved Theorem
\ref{param} for finite volume cone-manifolds with all cone angles $\leq 2\pi$ without the restriction on the tube radius (\cite{Hodgson:Kerckhoff:cone}). More recently they have announced that the parameterization holds with the tube radius condition we give here (\cite{Hodgson:Kerckhoff:harmonic,Hodgson:Kerckhoff:tube}). (When the cone angle is zero the tube radius is infinite and the result holds.) The parameterization was extended to geometrically finite cone-manifolds in \cite{Bromberg:rigidity}.

\medskip

Although it is not necessary we simplify our notation by assuming that
all the cone angles are equal to a single cone angle $\alpha$. By Theorem
\ref{param} there is a neighborhood $V$ of $\Psi(g)$ where $\Psi$ is
invertible. Let $[g_t] = \Psi^{-1}(V \cap (t, \dots, t, X))$ with $[g] =
[g_\alpha]$.

We set notation that we will use throughout the remainder of the paper. Let $M_t = (M, g_t)$ be the one-parameter family of hyperbolic cone-manifolds coming from the cone-metrics $[g_t]$. Although the conformal boundary is a fixed conformal structure $X$, the projective boundary will change. Let $\Sigma_t$ denote the projective boundary of $M_t$. We label the connected components of the conformal boundary $X^1, \dots, X^k$ and the corresponding components of the projective boundary $\Sigma^1_t, \dots, \Sigma^k_t$.

If $\gamma$ is a simple closed curve in $M$ then $\cL_\gamma(g) = L_\gamma(g) + \imath \Theta_\gamma(g)$ is the complex length of the geodesic representative (if it exists) of $\gamma$ in $(M,g)$ where $L_\gamma(g)$ is the length of $\gamma$ and $\Theta_\gamma(g)$ is the twisting. Note that $\Theta_\gamma(g)$ is defined modulo the cone angle if $\gamma$ is a component of the cone singularity and modulo $2\pi$ if $\gamma$ is a smooth geodesic. For the one parameter family of metrics $g_t$ we simplify notation and write the complex length $\cL_\gamma(t)$. We also simplify notation by setting $L_{\cC}(t) = \underset{c \in \cC}{\Sigma} L_c(t)$.

We now state our first main result:
\begin{theorem}
\label{downtozero}
Let $M_\alpha \in GF(N, \cC)$ be a geometrically finite hyperbolic
cone-metric with cone angle $\alpha$. Suppose the tube radius of the cone singularity is $\geq
\sinh^{-1} \sqrt{2}$. Then there exists an
$\ell_0$ depending only on $\alpha$ such that if $L_c(\alpha) \leq \ell_0$ for all $c \in \cC$ then the one parameter family of cone-manifolds $M_t \in GF(N, \cC)$ is defined for all $t \leq \alpha$.
\end{theorem}

{\bf Remark.} The lower bound on the tube radius is not essential. If we eliminate the lower bound on the tube radius the theorem still holds although $\ell_0$ will then depend on both the cone angle and the tube radius of the cone singularity. We have chosen a fixed tube radius to simplify the exposition.

\medskip

In our next result we control the geometry of the geometrically finite ends as the cone angle decreases. In particular we will measure the distance between the projective boundaries of $M_\alpha$ and $M_t$. This distance is defined in the next section. We also note that $\|\Sigma^i_\alpha\|_\infty$ is the distance between $\Sigma^i_\alpha$ and the unique Fuchsian projective structure with conformal structure $X^i$. This is also defined in the next section.

\begin{theorem}
\label{bigschwarzbound}
There exists a $C$ depending only on $\alpha$, the injectivity radius of the unique hyperbolic metric on $X^i$ and $\|\Sigma^i_\alpha\|_\infty$ such that
$$d(\Sigma^i_\alpha, \Sigma^i_t) \leq CL_{\cC}(\alpha)$$
for all $t \leq \alpha$.
\end{theorem}

We can also control the complex length of geodesics in $M_t$.
\begin{theorem}
\label{controllengths}
For each $L>0$ there exists an $\epsilon >0$ and $A > 0$ such that if $\gamma$ is a simple closed curve in $M$ with $L_\gamma(\alpha) \leq L$ and $L_c(\alpha) \leq \epsilon$ for all $c \in \cC$ then
$$e^{-AL_{\cC}(\alpha)} L_\gamma(\alpha) \leq L_\gamma(t) \leq e^{AL_{\cC}(\alpha)} L_\gamma(\alpha)$$
and
$$(1-AL_{\cC}(\alpha))\Theta_\gamma(\alpha) \leq \Theta_\gamma(t) \leq (1+AL_{\cC}(\alpha))\Theta_\gamma(\alpha)$$
for all $t \leq \alpha$.
\end{theorem}

We note that Theorems \ref{bigschwarzbound} and \ref{controllengths} are proved before Theorem \ref{downtozero} at least for all $t \leq \alpha$ where $M_t$ is known to exist. In fact these results are used in the proof of Theorem \ref{downtozero}.

\section{Deforming $(PSL_2\cx, X)$-structures}
Although there is a very general theory of $(G, X)$-structures, for
simplicity we will restrict to $(PSL_2\cx, X)$-structures here. In
fact for our purposes $X$ will either be $\chat$ or
$\hthree$. Then a $(PSL_2\cx, X)$-structure on a manifold is an atlas of
charts to $X$ with transition maps in $PSL_2\cx$; i.e. either a projective or
hyperbolic structure. We use the $(PSL_2\cx, X)$ notation simply so that we
can develop together the common elements of the deformation theory of
projective and hyperbolic structures.

An equivalent way to define a $(PSL_2\cx, X)$-structure is through a developing map and holonomy representation. A {\em developing map} $D$ is a local diffeomorphism from the universal cover $\tilde{M}$ to $X$ that commutes with a holonomy representation $\rho: \pi_1(M) \longrightarrow PSL_2\cx$. That is
\begin{equation}
\label{devmap}
D(\gamma(x)) = \rho(\gamma)(D(x))
\end{equation}
for all $\gamma \in \pi_1(M)$ and $x \in \tilde{M}$.

Let $\cD(M)$ be the space of developing maps for $M$ which we topologize with the compact-$C^\infty$ topology. We also define an equivalence relation for developing maps. We say $D_1 \sim D_2$ if there exists a diffeomorphism $\psi: M \longrightarrow M$ isotopic to the identity and element $\alpha \in PSL_2\cx$ such that $D_1 = \alpha \circ D_2 \circ \tilde{\psi}$. Let $D(M)$ be the quotient space $\cD(M)/\sim$.

To study one-parameter families of $(PSL_2\cx, X)$-structures we
need to make a definition about vector fields on $X$ and $M$. We say a
vector field $v$ on $X$ is {\em geometric} if the 
homeomorphisms in the one-parameter flow it defines are elements of
$PSL_2\cx$. As is well known the space of geometric vector fields is
the Lie algebra $sl_2\cx$. Geometric vector fields are analytic in the sense that any geometric vector field is determined uniquely by its germ at a single point. If $M$ has a $(PSL_2\cx, X)$-structure then a vector
field $v$ on $M$ is geometric if for every chart $\phi$, $\phi_* v$
is geometric.

A one-parameter family $M_t$ of $(PSL_2\cx,X)$-structures on $M$ can
be defined through a one parameter family of developing maps $D_t$
and holonomy representations $\rho_t$. By 
taking the derivative of $D_t$ we can define a family of vector
fields $v_t$ on the universal covers $\tilde{M}_t$. More precisely,
if $x$ is a point in $\tilde{M}$ then $D_t(x)$ is a smooth path in
$X$. The derivative $D'_t(x)$ will be a tangent vector
to the path at $D_t(x)$. We pull back this tangent vector to $T_x
\tilde{M}$ via $D_t$ to define the vector field $v_t$ at $x$.

Although these vector fields are defined on the differentiable manifold
$\tilde{M}$ the vector field $v_t$ has a special automorphic property in
the $(PSL_2\cx, X)$-structure on $M_t$. Explicitly, by differentiating
(\ref{devmap}) we see that for each $\gamma \in \pi_1(M)$ the vector field
$v - \gamma_* v$ is a geometric on $\tilde{M}_t$. We say that any
vector field that satisfies this relationship is {\em automorphic}.
To see how an automorphic vector field
describes the infinitesimal change in geometry we need to
discuss projective structures and hyperbolic structures separately.

\subsection{Projective structures}
A projective structure on a surface $S$ is a $(PSL_2\cx, \chat)$-structure.
As we noted when we first defined projective structures a projective structure also defines a conformal structure on $S$ and a fixed
conformal structure $X$ will have many distinct projective
structures. We let $P(X)$
denote the space of projective structures on $S$ with conformal
structure $X$. $P(X)$ inherits a topology as a subspace of the space of developing maps $D(S)$.  We will only be interested in projective structure
deformations contained in $P(X)$. This greatly simplifies the theory.

The objects that distinguish different projective structures in $P(X)$
are holomorphic quadratic differentials. In a local, conformal chart for
$X$ a holomorphic quadratic differential $\Phi$ has
the form $\Phi(z) = \phi(z) dz^2$ where $\phi$ is a holomorphic
function. On this local chart a conformal metric is determined by a
positive, real valued function $\sigma$ and has the form $\sigma(z)^2 dz
d\zbar$. On a Riemann surface there is at most one complete, conformal hyperbolic metric and this will always be the metric we will use. In particular if $\sigma$ is the hyperbolic metric on $X$ then $\|\Phi(z)\|_X = \sigma^{-2} |\phi(z)|$ is a well defined function on $S$ which we define to be the point-wise norm of $\Phi$ with respect to the
$\sigma$ metric. When it is clear which conformal structure is determining the metric we will drop the subscript and write the norm $\|\Phi(z)\|$.

Our first construction of a holomorphic quadratic differential will come from a conformal vector field $v$ on a projective structure $\Sigma$.
In a local chart, $v$ has the form $v(z) =
f(z) \ddz$. Since $v$ is conformal, $f$ is a holomorphic function. The {\em Schwarzian derivative} $Sv$ of
$v$ is a quadratic differential defined in a local projective chart by
$$Sv(z) = f_{zzz}(z) dz^2.$$
Note that this will only be a well defined quadratic differential if
the derivative of $f$ is taken in projective charts. If it is taken in
an arbitrary conformal chart the equation will not define a quadratic
differential. For projective structures a geometric vector field is usually called {\em projective}. As is well known $v$ will be projective if and only if $f(z)$ is a quadratic polynomial in $z$. In particular, if $Sv \equiv 0$ if and only if $v$ is a
projective vector field. The flow of a projective field preserves the projective structure so the Schwarzian measures the infinitesimal change in
projective structure. 

Note that there are no global conformal vector fields on a closed Riemann surface of genus $> 1$. The conformal vector fields we will be interested in are automorphic vector fields $v$ on the universal cover $\tilde{\Sigma}$ of our projective structure $\Sigma$. Then $Sv$ will be a holomorphic quadratic differential on $\tilde{\Sigma}$. However, by the automorphic property $Sv$ will be equivariant and descend to a quadratic differential on $\Sigma$.

The second holomorphic quadratic differential we will construct will measure the distance between two projective structures $\Sigma_0$ and $\Sigma_1$ in $P(X)$. There
is an obvious map $f$ between $\Sigma_0$ and $\Sigma_1$, namely the unique conformal
map. It is the existence of this map that simplifies the deformation theory
of projective structures in $P(X)$. The Schwarzian
derivative of $f$, defined using projective charts for $\Sigma_0$ and
$\Sigma_1$, is the quadratic differential 
\begin{equation}
\label{schwarziandef}
Sf = \left[\left(\frac{f_{zz}}{f_z}\right)_z -
\frac12\left(\frac{f_{zz}}{f_z}\right)^2\right]dz^2
\end{equation}
Again we must use projective charts for this equation to give a well
defined quadratic differential.

Conversely given any projective structure $\Sigma_0 \in P(X)$ and 
any holomorphic quadratic differential $\Phi$ on $X$, there exists a
unique projective structure $\Sigma_1$ in $P(X)$ such that 
$\Phi = Sf$. In particular, after fixing $\Sigma_0$ as a basepoint
there is a canonical isomorphism from $P(X)$ to  $Q(X)$, the space of
holomorphic quadratic differentials on $X$. The space $Q(X)$ is a
finite dimensional vector space. The $L^\infty$-norm on quadratic
differentials makes $Q(X)$ a normed vector space. The identification
of $P(X)$ with $Q(X)$ gives $P(X)$ a euclidean metric which will not
depend on the choice of basepoint. In this metric the distance
between two projective structures $\Sigma_0$ and $\Sigma_1$ will be
$d(\Sigma_0, \Sigma_1) = \|Sf\|_\infty$.

A projective structure is {\em Fuchsian} if its developing map is a
homeomorphism onto a round disk in $\chat$. The uniformization theorem
implies that there is a unique Fuchsian projective structure in $P(X)$.
As we  will frequently need to know the distance between an arbitrary
projective structure $\Sigma \in P(X)$ and the unique Fuchsian
element $\Sigma_F$ we define $\|\Sigma\|_\infty = d(\Sigma, \Sigma_F)$.

A vector space is its own tangent space so the derivative of a smooth path $\Sigma_t$ in $P(X)=Q(X)$ will give also be smooth path of holomorphic quadratic differentials $\Phi_t$ on $X$. To see this more explicity we let $\tilde{X}$ be the universal cover of the conformal
structure $X$ and choose conformal developing maps $D_t: \tilde{X} \longrightarrow \chat$. The vector
fields $v_t$ obtained by differentiating $D_t$ will be conformal and automorphic on $\tilde{\Sigma}_t$ therefore $Sv_t$ will be a holomorphic quadratic differential on $X$. By noting that on suitably chosen local charts $f_t=D_t \circ D^{-1}_0$ and differentiating (\ref{schwarziandef}) we see that $Sv_t$ is exactly $\Phi_t$. This implies the following proposition:

\begin{prop}
\label{schwarzpathbound}
The length of a smooth path, $\Sigma_t$ with $a< t < b$, in $P(X)$ is
$$\int_a^b \|\Phi_t\|_\infty dt.$$
\end{prop}

\subsection{Hyperbolic structures}
\label{Hyperbolic structures}
For a family of hyperbolic structures there is no obvious choice of
maps between the structures that plays the role of the conformal map
in the deformation theory of projective structures. However, recent
work has show that there is a canonical choice of an automorphic
vector field describing an infinitesimal deformation of the hyperbolic
structure. We describe this in more detail below. Our review will be brief. See \cite{Hodgson:Kerckhoff:cone} for more details.

Assume $M$ has a fixed hyperbolic metric coming from a developing map
$D$ with holonomy $\rho$. In hyperbolic space a geometric vector field is an infinitesimal isometry or {\em Killing field}. Let $E$ and $\tilde{E}$ be the bundles
over $M$ and $\tilde{M}$, respectively, of germs of Killing fields. For 
$\tilde{M}$ the developing map identifies germs of 
Killing fields at a point in $\tilde{M}$ with a Killing field on
$\hthree$ so  $\tilde{E}$ has a global product structure,
i.e. $\tilde{E} \isom \tilde{M} \times sl_2\cx$. Then $E$ is the
quotient of $\tilde{E}$ by the action of $\pi_1(M)$ where the action
on the first factor is by 
deck transformations and on the second by the holonomy representation.
The product structure on $\tilde{E}$ defines a flat connection which
descends to the quotient $E$. This flat connection has a covariant
derivative $d$ which we use to define the deRham cohomology groups
$H^i(M; E)$. As we shall see next, an automorphic vector field on
$\tilde{M}$ determines a cohomology class in $H^1(M; E)$. 

Given a vector field $v$ on $M$ (or $\tilde{M}$) the {\em canonical
lift} $s$ of $v$ is the section of $E$ (or $\tilde{E}$) determined by
the relationship that $s(p)$ is the unique Killing field that agrees
with $v$ at $p$ and whose curl agrees with the curl of $v$ at $p$. If a section $s$ of $\tilde{E}$ is the canonical lift of an automorphic vector field on $\tilde{M}$ then $s - \gamma_* s$ will be a constant section. 
Therefore $\omega = ds$ will descend to a
closed $E$-valued 1-form on $M$ which determines a cohomology class in $H^1(M;E)$.

A smooth path of developing maps $D_t$ in $\cD(M)$ defines an automorphic vector field in the following way. Let $x$ be a point in $\tilde{M}$. Then $D_t(x)$ is a smooth path in $\hthree$. Let $v_t(x)$ be the pull back of the tangent vector at the point $D_t(x)$ of this path by the developing map $D_t$. This defines a one parameter family of vector fields on $\tilde{M}$. On the hyperbolic structure $M_t$ defined by the developing map $D_t$ the vector field $v_t$ will be automorphic.

If we let $s_t$ be the canonical lift of $v_t$ then $\omega_t = ds_t$ is a family of $E_t$-valued 1-forms. Furthermore if the paths $D_t$ and $D'_t$ are equivalent in $D(M)$ the $\omega_t$ and $\omega'_t$ will be cohomologous. Therefore the derivative of a path in $D(M)$ is a one parameter family of cohomology classes in $H^1(M_t; E_t)$. This cohomology class plays the role of the holomorphic quadratic differential in the study of projective structures.

For many calculations we will take advantage of the complex structure
on $E$. In particular, the Lie algebra $sl_2\cx$ has a complex
structure that can be interpreted geometrically. If $v$ is a Killing
field on $\hthree$ then $\curl v$ is also a Killing field and $\curl
\curl v =  -v$. Therefore taking the curl of $v$ is equivalent to
multiplying by $\imath$. For a section $s$ of $E$ this leads us to
define $\imath s$ by the relationship $\imath s(p) = \curl
(s(p))$. We make a similar definition for $E$-valued $n$-forms. A section $s$ of $E$ is {\em real} if the Killing field
$\imath s(p)$ is zero at $p$, while $s$ is {\em imaginary} if $s(p)$
is zero at $p$. Every section $s$ has a unique decomposition into a
real section $\Re s$ and an imaginary section $\Im s$. A real section determines a vector field by the formula $v(p) =
(s(p))(p)$ and vice versa. If $s$ is an imaginary section then $\imath s$ is a real section so the formula $v(p) = (\imath s(p))(p)$ also identifies each imaginary section with a vector field. Returning to a general section $s$ we have a map $s \mapsto (\Re s, \Im s)$. If we view both $\Re s$ and $\Im s$ as vector fields this defines an isomorphism $E \longrightarrow TM \oplus TM$. The canonical lift of a vector field $v$ is then $(v, -\curl v)$ under this isomorphism.

This identification of $E$ with $TM \oplus TM$ gives $E$ a natural metric; we simply use the hyperbolic metric on each copy of $TM$. This metric defines an isomorphism from $E$ to the dual bundle $E^*$. For an $E$-valued $k$-form $\alpha$ we let $\alpha^\sharp$ be the image of $\alpha$ in $E^*$ under this isomorphism while if $\alpha$ is an $E^*$-valued $k$-form we let $\alpha^\flat$ be the image of $\alpha$ under the inverse of the isomorphism. The bundle $E^*$ has a flat connection $d^*$ and we define $\del \alpha = (d^* \alpha^\sharp)^\flat$. The formal adjoint for $d$ defined on $k$-forms is $\delta = (-1)^k * \del *$ where $*$ is the Hodge star operator. We also define the Laplacian $\Delta = d\delta + \delta d$.

In \S1 of \cite{Hodgson:Kerckhoff:cone} there are explicit formulas for $d$ and $\delta$ in local coordinates in terms of the Riemannian connection $\nabla$ and algebraic operators. Let $\{e_i\}$ be an orthonormal frame field with dual co-frame field $\{\omega_i\}$. Then
\begin{equation}
\label{dformula}
d = \underset{i}{\Sigma} \omega^i \wedge (\nabla_{e_i} + \ad(E_i))
\end{equation}
and
\begin{equation}
\label{deltaformula}
\delta = \underset{j}{\Sigma} i(e_j)(\nabla_{e_j} - \ad(E_j)).
\end{equation}
Here $i()$ is the interior product on forms. The operator $\ad(E_i)$ takes a Killing field field $Y$ to the Killing field $[E_i, Y]$. We also decompose $d$ and $\delta$ into the real and imaginary parts. Namely we let $D = \Re d$, $T= \Im d$, $D^* = \Re \delta$ and $T^* = \Im \delta$. Formulas for $D$, $T$, $D^*$ and $T^*$ follow easily from \eqref{dformula} and \eqref{deltaformula} since $\nabla_{e_i}$ is a real operator and $\ad(E_i)$ is an imaginary operator. In particular $T$ and $T^*$ are algebraic operators and therefore easy to calculate. It is also worth noting that $\del = D - T$. That is the flat connection on $E^*$ is the "conjugate" of the flat connection on $E$.

We make three more definitions that will be useful later. Let
$$\|\alpha\|^2 = \alpha \wedge *\alpha^\sharp.$$
Strictly speaking $\|\alpha\|^2$ is a real valued 3-form while $*(\alpha \wedge *\alpha^\sharp)$ is a function. We will abuse notation and use $\|\alpha\|^2$ to refer to both the 3-form and the function. It should be clear from context which meaning is correct. A vector 
field $v$ is {\em harmonic} if $\Delta s = 0$ where $s$ is the
canonical lift of $v$. A closed $E$-valued 1-form is a {\em Hodge form} if $\omega = ds$ where $s$ is the canonical lift of an harmonic, automorphic,
divergence free vector field on $\tilde{M}$. By the work in \S2 of \cite{Hodgson:Kerckhoff:cone} $\omega$ is a Hodge form if and only if $\omega$ is closed and co-closed and the real and imaginary parts of $\omega$ are symmetric and traceless vector valued 1-forms.

There is a nice formula for the $L^2$-norm of a Hodge form:
\begin{theorem}[\cite{Hodgson:Kerckhoff:cone}]
\label{normform}
Let $\omega$ be an $E$-valued Hodge form on a compact hyperbolic
3-manifold $M$ with boundary. Then
$$2\int_M \|\omega\|^2  = \int_{\del M} \imath \omega \wedge
\omega^\sharp$$
where $\del M$ is oriented with inward pointing normal.
\end{theorem}

\subsection{Extending deformations to the projective boundary}

Let $M$ be a hyperbolic 3-manifold with projective boundary
$\Sigma$. Together $M$ and $\Sigma$ form a differentiable 3-manifold
with boundary so if $v$ is a vector field on $M$ we can discuss its
continuous extension to $\Sigma$ and vice versa. We will always want
the extended vector field to be tangent to the boundary. We will use
this notion to discuss extending $E$-valued 1-forms on $M$ to
holomorphic quadratic differentials on $\Sigma$. Essentially an
$E$-valued 1-form extends to continuously to a holomorphic quadratic
differential if the vector field that generates the 1-form extends
continuously to a vector field that generates the quadratic differential.

Our vector fields will in general be automorphic vector fields on the
universal cover. However, the extension property is a local one so we
will work in an open neighborhood $V$ contained in $\hthree$ and with
boundary $\del V$ contained in $\chat$. Let $\omega$ be an
$E$-valued 1-form on $V$ and $\Phi$ a holomorphic quadratic
differential on $\del V$. Then $\omega = ds$ where $s$ is an $E$-valued section on $V$ and $\Phi = Sv_\infty$
where $v_\infty$ is a conformal vector field on $\del V$. Neither $s$
nor $v_\infty$ are uniquely determined; we can add a constant section to $s$ and a projective vector field to $v_\infty$. We say that $\omega$
extends continuously to $\Phi$ if $s$ and $v_\infty$ can be chosen
such that $\Re s$ extends continuously to $v_\infty$ and $-\Im s$ extends
continuously to $\imath v_\infty$.

Returning to our hyperbolic 3-manifold $M$ with projective boundary
$\Sigma$, the $E$-valued 1-form $\omega$ extends continuously to the
holomorphic quadratic differential $\Phi$ if it does so in a
neighborhood of every point on $\Sigma$. In general, an $E$-valued
1-form is {\em conformal at infinity} if it is cohomologous to a 1-form that
extends continuously to a holomorphic quadratic differential on
$\Sigma$.

Recall the one-parameter family of cone-manifolds $M_t$ that we defined in the previous section. The derivative of this path will be a path of cohomology classes $[\omega_t]$ in $H^1(M_t; E_t)$. The derivative of the projective boundary will be a path of holomorphic quadratic differentials $\Phi_t$ in $P(X)$. The following Hodge theorem is Theorem 4.4 in \cite{Bromberg:rigidity} plus Theorem \ref{context} in the appendix of this paper.

\begin{theorem}
\label{hodge}
The cohomology class $[\omega_t]$ is represented by a Hodge form $\omega_t$ which extends continuously to $\Phi_t$ on $\Sigma_t$.
\end{theorem}

Theorem \ref{normform} tells us how to calculate the $L^2$-norm of a Hodge form on a compact hyperbolic 3-manifold with boundary. We will need to calculate the $L^2$-norm of conformal Hodge forms on geometrically finite manifolds.

\begin{theorem}
\label{conformalnormform}
Let $M$ be a hyperbolic 3-manifold with boundary such that the union of $M$ with its projective boundary is compact. If $\omega$ is a conformal Hodge form that extends continuously to a holomorphic quadratric differential on the projective boundary then
$$2\int_M \|\omega\|^2 = \int_{\del M} \imath \omega \wedge \omega^\sharp$$
where $\del M$ is oriented with inward pointing normal.
\end{theorem}

\section{Tubes and half spaces}
\label{tubehalf}
In this section we make a digression from studying families of cone-manifolds to prove some results about a single hyperbolic cone-manifold $M_\alpha = (M, g)$ with cone-metric $g$ and all cone angles $\alpha$. Our goal is to find a constant $\ell_0$ such that if the length of the cone singularity is less than $\ell_0$ it will have a "large" tubular neighborhood and this neighborhood will be disjoint from any embedded half space in the geometrically infinite ends. We will prove a succession of results, each producing its own constant. At the end of the section we will simply take the minimum of these constants to find a single constant which will be used throughout the rest of the paper.
 
We first review an estimate of Hodgson and Kerckhoff on the size of
embedded tubes about the cone locus. These should be thought of as cone-manifold versions of the Margulis lemma with explicit constants. The main difficultly is of course the cone singularity. As a first stab at proving these results one might hope to smooth the cone-metric and then apply the Margulis lemma for manifolds with pinched negative curvature. There are two problems with this approach. If the cone angles are $< 2\pi$ then it may not be possible to smooth the metric to a negatively curved metric. On the other hand if there a cone angles $> 2\pi$ we can always smooth the metric to one that is negatively curved however we cannot bound the amount of negative curvature required even if all the cone angles are bounded. To get around both of these problems we assume a priori that the cone singularity has a tubular neighborhood of definite size. In practice this is not much of a restrictions since in most applications the hyperbolic cone-manifolds arise from smooth hyperbolic structures where the standard Margulis lemma applies.

\begin{prop}
\label{margulis}
Let $M_\alpha = (M,g)$ be a hyperbolic cone-manifold with all cone angles
$\alpha$ and let $\gamma$ be a closed non-singular geodesic in
$M_\alpha$. Suppose the tube radius of the cone singularity is $\geq 
R$. Then there exists an
$\ell_1 > 0$ depending only on $\alpha$ and $R$ such that if $\gamma$ is a
closed geodesic with $L_\gamma(g) \leq \ell_1$ and $L_c(g) \leq \ell_1$ for all
$c \in \cC$ then $\gamma$ has an embedded tube of radius
$R$ which is disjoint from the 
$R$-neighborhood of the cone
singularity.
\end{prop}

{\bf Proof.} We can construct a complete Riemannian metric $h$ on $M$
such that $h$ agrees with $g$ outside the $R/2$-neighborhood of the cone singularity and such that $h$ has pinched
negative sectional curvature (see \cite{Kojima:cone} Theorem 1.2.1). If $L_\gamma(g) \leq R/2$ then $\gamma$ will be disjoint from the $R$ neighborhood of $\cC$ and therefore $\gamma$ will also be a geodesic in the $h$ metric with $L_\gamma(g) = L_\gamma(h)$. Furthermore there is a universal bound
on the sectional curvature of $h$ depending only on our choice of tube
radius $R/2$. 

Let $M^\epsilon_h$ be the $\epsilon$-thin part of $M$ for the $h$
metric. That is $M^\epsilon_h$ is the subset of $M$ consisting of
those points whose injectivity radius is $< \epsilon$. By the Margulis
lemma (\cite{Ballman:Gromov:Schroeder}) there is an $\epsilon_0$, depending only on the curvature
bounds, such that each component of $M^{\epsilon_0}_h$ has virtually
nilpotent fundamental group. Since $M$ is an orientable, hyperbolizable
3-manifold the only possible virtually nilpotent subgroups of $M$ are
$\integers$ and $\integers \oplus \integers$, and the second case will
only occur at peripheral tori. 

Let $c$ be a component of $\cC$ and $V_c$ the $R$-neighborhood of $c$ in the $g$ metric. Choose $\delta_1$ such that if
$L_c(g) \leq \delta_1$ then $V_c$ will be contained in the
$\epsilon_0$-thin part of the $g$ metric.  This is the one place where
our choice depends on $\alpha$ for as the cone angle increases,
$\delta_1$ must decrease. The $h$ metric will decrease the injectivity
radius in $V_c$ so $V_c$ will also be contained in $M^{\epsilon_0}_h$.

Next choose $\delta_2$ such that if $L_\gamma(g) \leq \delta_2$ then the
$R$-neighborhood $V_\gamma$ of $\gamma$ in the $h$ metric is also contained in
$M^{\epsilon_0}_h$. Note that $\delta_2$ will only depend on the
curvature bounds and not on the cone angle. Any component of
$M^{\epsilon_0}_h$ with fundamental group $\integers \oplus \integers$
does not contain any closed geodesics, therefore the component of
$M^{\epsilon_0}_h$ containing $\gamma$ must be a solid torus with
fundamental group $\integers$. This implies that $V_\gamma$ is a solid torus and is disjoint from $V_c$. Since $V_\gamma$ is disjoint from $V_c$, the $g$ and $h$ metrics agree on $V_\gamma$ so $\gamma$ has embedded tubular neighborhood of radius $R$ in the original metric if $L_\gamma(g) \leq \ell_1$ where $\ell_1 = \min \{R/2, \delta_1, \delta_2\}$. \qed{margulis} 

Combining this proposition with Theorem 4.4 in \cite{Hodgson:Kerckhoff:dehn} we have:

\begin{prop}
\label{firsttubebound}
Let $M_\alpha = (M,g)$ be a hyperbolic cone-manifold with all cone angles
$\alpha$ and $\gamma$ a closed geodesic in $M_\alpha$. Suppose the tube radius of the cone singularity is $\geq R$, $L_c(g) \leq \ell_1$ for each $c \in \cC$ and $L_\gamma(g) \leq \ell_1$. Then for each $c \in \cC$, $c$ and the geodesic $\gamma$ have disjoint tubular neighborhoods such that the area of the boundary tori is $\geq 1.6978\frac{\sinh^2 R}{\cosh (2R)}$.
\end{prop}


We next prove similar results for hyperbolic half spaces.

\begin{lemma}
\label{allhalfembed}
Let $D$ be an embedded round disk in the projective boundary of
$M_\alpha$. Then there is an embedded hyperbolic half space $H$ in $M$
whose projective boundary is $D$. 
\end{lemma}

{\bf Proof.}  In general a hyperbolic half space $H$ is bounded by a hyperbolic plane $P$ and it projective boundary, a round disk $D$. The half space is foliated by planes $P_d$ of constant curvature where $P_d$ is the set of points a distance $d$ from $P$. If $D$ is an embedded round disk on the boundary of a geometrically finite cone-manifold then the planes $P_d$ will be embedded for large $d$. Let $d'$ be the inf of all such $d$. If $d'>0$ then the metric closure of $\underset{d > d'}{\cup} P_d$ will have strictly concave boundary so $M'= M \backslash \underset{d>d'}{\cup} P_d$ will have strictly convex boundary $P_{d'}$. This implies that $P_{d'}$ is embedded hence we must have $d' =0$. The only way $P_0$ cannot be embedded is for it to intersect an element $c$ of the cone-singularity. In this case $c$ must be tangent to $P_0$ and therefore contained in $P_0$. Since $P_0$ does not contain a closed geodesic it cannot intersect $c$. Therefore $P_0$ and hence all of $H$ must be embedded. \qed{allhalfembed}

Next we see that these embedded half spaces do not intersect the tubular neighborhood of the cone singularity if it is sufficiently short.

\begin{prop}
\label{halfdisjointtube}
Let $M_\alpha= (M,g)$ be a hyperbolic cone-manifold with all cone angles
$\alpha$. Suppose that the tube radius of the cone-singularity is
$\geq R$. Then there exists
an $\ell_2$ depending only on $\alpha$ and $R$ such that if $L_c(g) \leq \ell_2$
for all $c \in \cC$ then any embedded half space $H$ is disjoint from a tubular neighborhood of $c$ with the area of the boundary torus $\geq 1.6978\frac{\sinh^2 R}{\cosh (2R)}$.
\end{prop}

{\bf Proof.}  We first show that $\ell_2$ can be chosen such that if $L_c(g) \leq \ell_2$ for all $c \in \cC$ then $H$ does not intersect the $R$-neighborhood of the cone-singularity.  We need the following simple geometric fact. If $U$ is a tube of radius $R'>R$ and $H$ intersects the $R$-neighborhood of the core curve of the tube then the area of the intersection of $H$ with $\del U$ is bounded below by a function $A(R')$ with $A(R') \rightarrow \infty$ as $R' \rightarrow \infty$. Note that $\area(\del U) > A(R')$.

Choose $R'$ such that $A(R') = 1.6978\frac{\sinh R}{\cosh(2R)}$ and choose $\delta$ such that
$$\alpha \delta \sinh R' \cosh R' = 1.6978\frac{\sinh R}{\cosh(2R)}.$$
Then let $\ell_2 = \min \{\ell_1, \delta\}$. If $L_c(g) \leq \ell_2$ for all $c \in \cC$ then by Proposition \ref{firsttubebound} the $R'$-neighborhood of each $c \in \cC$ will be an embedded tube $U$. Furthermore
$$\area(\del U) = \alpha L_c(g) \sinh R' \cosh R' \leq 1.6978\frac{\sinh R}{\cosh(2R)}.$$
If $H$ intersect the $R$-neighborhood of $c$ then
$$\area(\del U)  > A(R') = 1.6978\frac{\sinh R}{\cosh(2R)}.$$
This contradiction implies that $H$ does not intersect the $R$-neighborhood.  

The rest of the proof is similar to the proof of Theorem 4.4 in
 \cite{Hodgson:Kerckhoff:dehn} and we will only sketch it. Define $R_m$ to be
 the maximal radius such that the tube $U_m$ of radius $R_m$ about $c$ is embedded and disjoint from $H$. We can assume that $T=\del U_m$ intersects the hyperbolic plane $P=\del H$ for otherwise we can simply apply Proposition \ref{firsttubebound}. Except for possible self-tangencies $T$ will be embedded and $P$ will be tangent to $T$ at a point $p$. Let $B$ be the ball of radius $R$ contained in $H$ and tangent to $T$ at $p$.

We now lift to the universal cover. Since $M$ is hyperbolizable, any component (again ignoring self-tangencies) $\tilde{T}$ of the pre-image of $T$ in $\tilde{M}$ will be a topological plane and the stabilizer of $\tilde{T}$ will be a $\integers \oplus \integers$ subgroup $\Gamma_T$ of $\pi_1(M)$. Let $\tilde{p}$ be a point in the pre-image of $p$ contained in $\tilde{T}$ and let $\tilde{B}$ be the component of the pre-image $B$ that is tangent to $\tilde{T}$ at $\tilde{p}$. Let $C$ be the orthogonal projection of $\tilde{B}$ onto $\tilde{T}$. One needs to be careful here to make sure that the cone singularity does not interfere with this orthogonal projection. It is at this point that we refer to Hodgson and Kerckhoff. In particular they show that $C$ is well defined and disjoint from its translates under the action of $\Gamma_T$. This implies that $\area(T) \geq \area(C)$. Hodgson and Kerckhoff also show $\area(C) \geq 1.6978\frac{\sinh^2 R}{\cosh (2R)}$. Therefore $\area(T) \geq 1.6978\frac{\sinh^2 R}{\cosh (2R)} $ as desired. \qed{halfdisjointtube}

We now define two constants determined by a cone metric $g$. For $c \in \cC$ let $R^c_g$ be chosen such that
$$\alpha L_c(g) \sinh R^c_g \cosh R^c_g = \frac{1}{2}.$$
Similarly if $\gamma$ is a closed geodesic let $R^\gamma_g$ be chosen such that
$$2\pi L_\gamma(g) \sinh R^\gamma_g \cosh R^\gamma_g = \frac{1}{2}.$$
Let $U^c_g$ and $U^\gamma_g$ be the $R^c_g$ and $R^\gamma_g$ neighborhoods of $c$ and $\gamma$, respectively, and let $U^{\cC}_g = \underset{c \in \cC}{\cup} U^c_t$. Note that the area of both $\del U^c_g$ and $\del U^\gamma_g$ is $\frac{1}{2}$.

In our next result we summarize the work of this section for a fixed choice of minimal tube radius. Our choice, although essentially arbitrary, will simplify some of the constants in the rest of the paper.

\begin{theorem}
\label{tubebound}
Let $M_\alpha = (M,g)$ be a hyperbolic cone-manifold with all cone angles $\alpha$ and assume that the tube radius of the cone singularity is greater than $\sinh^{-1} \sqrt{2}$. Then there exists an $\ell_0$ depending only on $\alpha$ such that if $L_c(g) \leq \ell_0$ for all $c \in \cC$ then the $U^c_g$ are embedded tubular neighborhoods, each pairwise disjoint and disjoint from any embedded half space and $R^c_g \geq \sinh^{-1} \sqrt{2}$. Furthermore if $\gamma$ is a closed geodesic with $L_\gamma(g) \leq \ell_0$ then $U^\gamma_g$ is also an embedded tubular neighborhood disjoint from the $U^c_g$. 
\end{theorem}

{\bf Proof.} Let $\ell_1$ and $\ell_2$ be the constants given by Propositions \ref{firsttubebound} and \ref{halfdisjointtube}, respectively, with $R = \sinh^{-1} \sqrt{2}$. Choose $\ell_3$ such that 
$$\alpha \ell_3 \sinh R \cosh R = \alpha \ell_3 (\sqrt{2}) (\sqrt{3}) = \frac{1}{2}.$$
Then $\ell_0 = \min \{\ell_1, \ell_2, \ell_3\}$ is the desired constant.
\qed{tubebound}

For each point $p$ in the projective boundary we will also need to
estimate the size of the largest embedded round disk containing $p$. Here
size will be measured by comparing the hyperbolic metric on the round
disk to the hyperbolic metric on the entire projective boundary. By
the Schwarz lemma the metric will always be bigger on the disk so we
want to find a disk where we can bound the ratio of the two
metrics. This bound will depend both on the injectivity radius of the
hyperbolic metric and on the deviation of the projective boundary from
a Fuchsian projective structure.

We begin with a simple lemma about hyperbolic geometry.

\begin{lemma}
\label{translate}
Let $D$ be a round disk in $\chat$ with hyperbolic metric $\sigma$.
Let $\gamma$ be an isometry of ${\mathbb H}^2$ with translation length $\geq \ell$. For every $z \in D$ there exists a round disk $D' \subset D$ such that $D' \cap \gamma(D') = \emptyset$ and $\sigma'(z) = \coth(\ell/4) \sigma(z)$, where $\sigma'$ is the hyperbolic metric for $D'$.
\end{lemma}

{\bf Proof.} Without loss of generality we assume that $D = \Delta$, the unit disk in $\cx$, and $z = 0$. Then let $D'$ be the euclidean disk of radius $\tanh (\ell/4)$ centered at $0$. The hyperbolic diameter of $D'$ in the $\sigma$ metric is $\ell/2$ so $D' \cap \gamma(D') =\emptyset$. Finally $\sigma(z) = \frac{2}{1 - |z|^2}$ while $\sigma'(z) = \frac{2\tanh(\ell/4)}{\tanh^2(\ell/4) - |z|^2}$ so
$$\sigma'(0) = \coth(\ell/4) \sigma(0)$$
as desired. \qed{translate}

Next, we use the previous lemma to estimate the size of embedded round
disks in a projective structure.

\begin{prop}
\label{embeddisks}
Let $\Sigma$ be a projective structure, $\sigma$ the hyperbolic metric
on $\Sigma$ and $\kappa$ the injectivity radius of $\sigma$. Then every $z \in \Sigma$ is contained in an embedded round disk $D$ in $\Sigma$ such that 
$$\sigma_D(z) < \sigma(z) \coth (\kappa/2)
\sqrt{1+2\|\Sigma\|_\infty}$$
where $\sigma_D$ is the hyperbolic metric on $D$.
\end{prop}

{\bf Proof.} Let $X$ be the conformal structure induced by $\Sigma$
and let $\Sigma_F$ be the unique Fuchsian structure in $P(X)$. Then
$\tilde{\Sigma}_F$ is projectively isomorphic to $D$ and the group of
deck transformation $\Gamma$ satisfies the conditions of Lemma
\ref{translate}. In particular for each $z \in \tilde{\Sigma}_0$ there
is a round disk $D'$ containing $z$ such that for each $\gamma \in
\Gamma$, $D \cap \gamma(D) = \emptyset$ and 
\begin{equation}
\label{transeq}
\sigma_{D'}(z) = \coth(\kappa/2) \sigma(z)
\end{equation}
where $\sigma_{D'}$ is the hyperbolic metric on $D'$.
Therefore $D'$ descends to an embedded round disk in
$\Sigma_F$.

Let $f$ be the unique conformal map from $\Sigma_F$ to $\Sigma$. Then $f(D')$ will be
will be an embedded topological disk but not a round
disk in the projective structure $\Sigma$. However by Theorem 4.2 in \cite{Anderson:projective} there exists a round disk in $D
\subset f(D')$ in $\Sigma$ with $z \in D$ such that
$$\sigma_D(z) \leq \sigma_{D'}(z) \sqrt{1 + 2\|Sf\|_{D',\infty}}.$$
By the Schwarz Lemma $\sigma_{D'} > \sigma$ so
$\|Sf\|_{D',\infty}  < \|Sf\|_{\Sigma, \infty} =  \|\Sigma\|_\infty$. Combining the two inequalities
with (\ref{transeq}) gives
$$\sigma_D(z) < \sigma(z) \coth (\kappa/2) \sqrt{1 + 2\|\Sigma\|_\infty}.$$
\qed{embeddisks}

\section{Bounding the length of geodesics}
We can now return to investigating the one parameter family of hyperbolic cone-manifolds $M_t$. 
As we shall see, the estimates we derive are simple consequences of the work of Hodgson and Kerckhoff.

Recall that the derivative of the path $M_t$ is a cohomology class in $H^1(M_t; E_t)$ represented by a Hodge form $\omega_t$. Throughout this section we assume that the one-parameter family is defined for $t$ in a half-open interval $(\alpha', \alpha]$ and that at the starting structure $L_c(\alpha) \leq \ell_0$ and $R^c_\alpha > \sinh^{-1} \sqrt{2}$ for all $c \in \cC$.

In our first result we show that the tube radius does not decay and we bound the length of the cone singularity.

\begin{prop}
\label{conelength}
For all $c \in \cC$, $L_c(t) \leq L_c(\alpha) \leq \ell_0$, $R^c_t > \sinh^{-1} \sqrt{2}$ and
\begin{equation}
\label{singlength}
\frac{tL_c(\alpha)}{\alpha + 2L_c(\alpha)(\alpha^2 - t^2)} \leq
L_c(t) \leq \frac{t L_c(\alpha) }{\alpha - 2L_c(\alpha)(\alpha^2 - t^2)}
\end{equation}
if $t \in (\alpha', \alpha]$.
\end{prop}

{\bf Proof.} By assumption $L_c(\alpha) \leq \ell_0$ and $R^c_\alpha > \sinh^{-1} \sqrt{2}$. We will show that these two properties hold for all $t \in (\alpha', \alpha]$. If the tube radius condition does not hold, by the continuity of $R^c_t$, there exists a largest value $T < \alpha$ such that $R^c_T = \sinh^{-1}\sqrt{2}$. We will work by contradiction and show that such a $T$ cannot exist.

To do so we show that $L_c(t) \leq \ell_0$ for $t \in [T,\alpha]$, also working by contradiction. If this does not hold there exists a $T'$ with $L_c(T') = \ell_0$, $L'_c(T') < 0$ and $L_c(t) \leq \ell_0$ if $t \in [T',\alpha]$. Proposition \ref{tubebound} implies that the tubular neighborhood of $c$ of radius $R^c_t$ is embedded for $t \in [T, \alpha]$ and therefore 
by an estimate of Hodgson and Kerckhoff (Theorem 2.7 of \cite{Hodgson:Kerckhoff:dehn}) we have:
\begin{equation}
\label{blah}
\frac{L_c(t)}{t}\left(1 - \frac{1}{\sinh^2 R^c_t} \right) \leq
L'_c(t) \leq \frac{L_c(t)}{t} \left(1 + \frac{1}{\sinh^2 R^c_t}
\right).
\end{equation}
Since $\frac{1}{\sinh^2 R^c_t} \leq \frac{1}{2}$, the left hand side of (\ref{blah}) implies $L'_c(T') > 0$. This contradiction implies that $L_c(t) \leq \ell_0$ for all $t$ between $\alpha$ and $T$. Since $TL_c(T) \sinh R^c_T \cosh R^c_T = \frac{1}{2}$ our choice of $\ell_0$ implies that $R^c_T > \sinh^{-1} \sqrt{2}$. Again, this is a contradiction and hence $R^c_t > \sinh^{-1} \sqrt{2}$ for all $t \in(\alpha', \alpha]$. The previous argument also shows that $L'_c(t) > 0$ and therefore $L_c(t) \leq L_c(\alpha) \leq \ell_0$ for all $t \in (\alpha', \alpha]$.

Next we combine the inequality
$$\frac{1}{\sinh^2 R^c_t} \leq \frac{2}{\sinh R^c_t \cosh R^c_t}$$
and the equality
$$tL_t(c) \sinh R^c_t \cosh R^c_t = \frac{1}{2}$$
with (\ref{blah}) to obtain
\begin{equation}
\label{coneder}
\frac{L_c(t)}{t} \left(1 - 4tL_c(t) \right) \leq L'_c(t) \leq \frac{L_c(t)}{t} \left( 1 + 4t L_c(t) \right).
\end{equation}
To prove (\ref{singlength}) we need to integrate this inequality. To
do so we make the substitution $y(t) = \frac{L_c(t)}{t}$. Then the
first inequality of (\ref{coneder}) becomes: 
$$y + t\frac{dy}{dt} \geq y(1 - 4t^2 y).$$
Rearranging and integrating we get
\begin{eqnarray*}
\int_T^\alpha \frac{1}{y^2} \frac{dy}{dt} dt & \geq & -\int_T^\alpha 4tdt \\
-\frac{1}{y(\alpha)} + \frac{1}{y(T)} & \geq & \frac{4(T^2 - \alpha^2)}{2} \\
\frac{T}{L_c(T)} &\geq& 2(T^2 - \alpha^2) +
\frac{\alpha}{L_c(\alpha)}.
\end{eqnarray*}
This final inequality is equivalent to the second inequality of
(\ref{singlength}). The other inequality is derived
similarly. 
\qed{conelength}

The $L^2$-norm of $\omega_t$ will be infinite on all of $M_t$. However if we let $M'_t = M_t \backslash U^{\cC}_t$ then the $L^2$-norm will be bounded on $M'_t$.

\begin{prop}
\label{normbound}
$$\int_{M'_t} \|\omega_t\|^2 \leq L_{\cC}(t)^2$$
Futhermore for any $A$ and $R$ there exists a $K$ such that if for each $c \in \cC$, $V^c$ is a tubular neighborhood of $c$ with the $\area(\del V^c) \geq A$ and the radius of $V^c$ greater than $R$ then
$$\int_{M_t \backslash V^{\cC}} \|\omega_t\|^2 \leq K^2 L_{\cC}(t)^2$$
where $V^{\cC} = \underset{c \in \cC}{\cup} V^c$.
\end{prop}

{\bf Proof.} By Proposition \ref{conelength}, $R^c_t > \sinh^{-1} \sqrt{2}$ for all $c \in \cC$ and $t \in (\alpha', \alpha]$. By Proposition \ref{tubebound} the tube $U^c_t$ of radius $R^c_t$ about $c$ is embedded. Using Theorem \ref{conformalnormform} along with work of 
Hodgson and Kerckhoff (see (17) on p. 14 of \cite{Hodgson:Kerckhoff:dehn}) we see that
$$\int_{M'_t} \|\omega_t\|^2 \leq \underset{c \in
  \cC}{\Sigma}\frac{\cosh R^c_t \area (\del U^c_t)}{t^2 \sinh^3
  R^c_t(2\cosh^2 R^c_t + 1)}.$$
From the area formula for the tube boundary and our definition of $R^c_t$, $\area (\del U^c_t) =tL_t(c) \sinh
R^c_t \cosh R^c_t = \frac{1}{2}$. Using this and also the fact that $\sinh R^c_t > \sqrt{2}$ we have:
\begin{eqnarray*}
\frac{\cosh R^c_t \area(\del U^c_t)}{t^2 \sinh^3 R^c_t(2 \cosh^2
  R^c_t + 1)} & = & \frac{2L_c(t)^2 \cosh^3 R^c_t}{\sinh R^c_t (2\cosh^2
  R^c_t + 1)} \\
& \leq & L_c(t)^2
\end{eqnarray*}
and therefore
$$\int_{M'_t} \|\omega_t\|^2 \leq \underset{c \in \cC}{\Sigma} L_c(t)^2 \leq L_{\cC}(t)^2.$$
The proof of the more general inequality is essentially the same.
\qed{normbound}

{\bf Remark.} In \cite{Hodgson:Kerckhoff:dehn} instead of using the cone angle as the parameter for the family of the hyperbolic cone-manifolds they use the cone angle squared. This accounts for the difference in the constants in their paper and the constants in this paper.

\medskip

We can also bound the length of short curves that are not part of the
cone singularity:

\begin{prop}
\label{lengthbound}
If $\gamma$ is a simple closed curve in $M$ with
$L_\gamma(T) \leq e^{-4\alpha\ell_0}\ell_0$ for some $T \in (\alpha', \alpha]$, then
\begin{equation}
\label{injbound}
e^{-4\alpha L_{\cC}(\alpha)}L_\gamma(\alpha) \leq L_\gamma(t) \leq e^{4\alpha L_{\cC}(\alpha)}L_\gamma(\alpha)
\end{equation}
and
\begin{equation}
\label{thetabound}
\left(1-4\ell_0 L_{\cC}(\alpha)\right) \Theta_\gamma(\alpha) \leq \Theta_\gamma(t) \leq \left(1 + 4\ell_0 L_{\cC}(\alpha)\right) \Theta_\gamma(\alpha)
\end{equation}
for all $t \in (\alpha', \alpha]$.
\end{prop}

{\bf Proof.} By Proposition \ref{conelength}, $R^c_t > \sinh^{-1} \sqrt{2}$ for all $c \in \cC$ and $t \in (\alpha', \alpha]$. Therefore by Proposition \ref{tubebound}, if $L_\gamma(t) \leq \ell_0$ the tube $U^\gamma_t$ is embedded and contained in $M'_t$. 

We need to show that the $L^2$-norm of $\omega_t$ on $U^\gamma_t$ is bounded by the derivative $\cL'_\gamma(t)$. The essential estimates again come from \S2 of \cite{Hodgson:Kerckhoff:dehn}. In particular they show that on $U^\gamma_t$, $\omega_t$ can be decomposed as the sum of a certain Hodge form of standard type and a correction term. That is
$$\omega_t = \frac{\cL'_\gamma(t)}{2\cL_\gamma(t)} \omega_l + \omega_c$$
where $\omega_\ell$ is a radially symetric Hodge form. For an explicit description of $\omega_l$ see p. 9 of \cite{Hodgson:Kerckhoff:dehn}. For our purposes we only need the following two facts. First, by Lemma 2.5 in \cite{Hodgson:Kerckhoff:dehn}
$$\int_{U^\gamma_t} \|\omega_t\|^2 = \left(\frac{|\cL'_\gamma(t)|}{2L_\gamma(t)}\right)^2 \int_{U^\gamma_t} \| \omega_l\|^2 + \int_{U^\gamma_t} \|\omega_c\|^2$$
and second, by the formulas on p. 14 of \cite{Hodgson:Kerckhoff:dehn}
$$\int_{U^\gamma_t} \| \omega_l\|^2 =  \frac{\sinh R^\gamma_t}{\cosh R^\gamma_t}
\left(2 + \frac{1}{\cosh^2 R^\gamma_t}
\right) \area(\del U^\gamma_t).$$
Together this implies
$$ \int_{U^\gamma_t} \|\omega_t\|^2 \geq \frac{\sinh R^\gamma_t}{\cosh R^\gamma_t}
\left(2 + \frac{1}{\cosh^2 R^\gamma_t}
\right)\left(\frac{|\cL'_\gamma(t)|}{2L_\gamma(t)}\right)^2 \area(\del U^\gamma_t).$$
Since $\sinh R^\gamma_t > \sqrt{2}$ and $\area(\del U^\gamma_t) = \frac{1}{2}$, 
$$\frac{\sinh R^\gamma_t}{\cosh R^\gamma_t}\left(2 + \frac{1}{\cosh^2 R^\gamma_t}\right) \area(\del U^\gamma_t) \geq \frac{\sqrt{2}}{\sqrt{3}}\left(2 + \frac{1}{3}\right)\left(\frac{1}{2}\right) > \frac{1}{4}$$
We also know that
$$\int_{M'_t} \|\omega_t\|^2 \geq \int_{U^{\gamma}_t} \|\omega_t\|^2.$$
By Theorem \ref{normbound}
$$L_{\cC}(t)^2 \geq \int_{M'_t} \|\omega_t\|^2.$$
By combining these four inequalities we have
$$L_{\cC}(t)^2 \geq \left(\frac{|\cL'_\gamma(t)|}{4L_\gamma(t)} \right)^2$$
which rearranges to give
\begin{equation}
\label{derbound}
|\cL'_\gamma(t)| \leq 4L_\gamma(t) L_{\cC}(t)
\end{equation}
if $L_\gamma(t) \leq \ell_0$.

Next we show that if $L_\gamma(T) \leq e^{-4\alpha \ell_0} \ell_0$ then $L_\gamma(t) \leq \ell_0$ for all $t \in (\alpha', \alpha]$. Let $I$ be the largest interval containing $T$ such that $L_\gamma(t) \leq \ell_0$ if $t \in I$. We will show that $I$ is an open and closed subset of $(\alpha', \alpha]$. By the continuity of the $L_\gamma(t)$, $I$ is closed. Furthermore if $T'$ is the right endpoint of $I$ then either $L_\gamma(T') = \ell_0$ or $T' = \alpha$. 
Since $|L'_\gamma (t)| \leq |\cL'_\gamma(t)|$ and $L_{\cC} (t) \leq L_{\cC}(\alpha)$, (\ref{derbound}) becomes
\begin{equation}
\label{localineq}
|L'_\gamma(t)| \leq 4L_\gamma(t) L_{\cC}(\alpha)
\end{equation}
if $t \in I$. By integrating (\ref{localineq}) from $T$ to $T'$ we get
\begin{equation}
\label{globalineq}
e^{-4|T'-T|L_{\cC}(\alpha)}L_\gamma(T) \leq L_{\gamma}(T') \leq e^{4|T'-T|L_{\cC}(\alpha)}L_\gamma(T).
\end{equation}
Since $|T'-T| \leq \alpha$ with equality only holding if $T' = \alpha$ and $T=0$, the right hand side of (\ref{globalineq}) implies that either $L_\gamma(T') < \ell_0$, which contradicts the definition of $I$, or $T' = \alpha$. A similar argument shows that the left endpoint of $I$ is $\alpha'$ and $I = (\alpha', \alpha]$. Therefore we can integrate (\ref{localineq}) to get 
(\ref{injbound}).

To prove (\ref{thetabound}) we note that (\ref{derbound}) implies that
$$|\Theta'_\gamma(t)| \leq 4L_\gamma(t) L_{\cC}(t) \leq 4\ell_0 L_{\cC}(\alpha).$$
Integrating this inequality gives (\ref{thetabound}). \qed{lengthbound}

Recall that the {\em injectivity radius} at a point $x$ in a Riemannian manifold is the radius of the largest embedded ball centered at $x$ or equivalently half the length of the shortest geodesic arc with both endpoints at $x$. We define $\inj_x(t)$ to be the injectivity radius at $x$ for the $g_t$ metric on $M$.

\begin{prop}
\label{boundinj}
The injectivity radius $\inj_x(t)$ is bounded away from zero on $M'_t$ for all $t \leq \alpha$.
\end{prop}

{\bf Proof.} Note that we are measuring the injectivity radius using the metric $g_t$ on all of $M$ but only showing that it is bounded below on $M'_t$. Clearly as $t \rightarrow 0$ this will not be true by Proposition \ref{lengthbound}.

If the injectivity radius is not bounded below there are two possibilities.

First, there could be a simple closed curve $\gamma$ in $M$ such that $L_\gamma(t) \rightarrow 0$. This is not possible by (\ref{injbound}).

The second possibility is that there are points $x_t \in \del U^c_t$ for some $c \in \cC$ such that $\inj_{x_t}(t) \rightarrow 0$ as $t \rightarrow 0$. The tori $\del U^c_t$ have an induced Euclidean metric and it is not hard to see that the hyperbolic injectivity radius will decay to zero if and only if the Euclidean injectivity radius decays to zero.

The Euclidean metric on $\del U^c_t$ can be constructed by gluing together (possibly with a twist) the boundary components of a Euclidean cylinder of height $L_c(t) \cosh R^c_t$ and radius $t\sinh R^c_t$. The area of this cylinder $tL_c(t) \sinh R^c_t \cosh R^c_T$ is always $\frac{1}{2}$ so the injectivity radius will be bounded below if and only if the height $H(t)$ is bounded above and below.  By (\ref{singlength}), $\frac{L_c(t)}{t}$ is bounded above and below. Since $R^c_t \geq \sinh^{-1} \sqrt{2}$, $\frac{\cosh R^c_t}{\sinh R^c_t}$ is also bounded above and below and therefore so is $\frac{L_c(t) \cosh R^c_t}{t \sinh R^c_t}$. Finally by multiplying the numerator and denominator by $L_c(t) \cosh R^c_t$ we have that
$$ \frac{\left(L_c(t) \cosh R^c_t\right)^2}{tL_c(t) \sinh R^c_t \cosh R^c_t} = 4H(t)^2$$
is also bounded above and below, as desired.
\qed{boundinj}

We can also bound the length of arbitrary geodesics although the estimates are slightly different.

\begin{lemma}
\label{outofthin}
Let $(M,g)$ be a hyperbolic cone-manifold. For each $c \in \cC$ assume that $U^c_g$ is an embedded tubular neighborhood with $R^c_g > \sinh^{-1} \sqrt{2}$ and assume that $L_c(g) \leq \ell_0$. For each $L>0$ there exits an $\epsilon > 0$ such that if $\gamma$ is a closed geodesic in $M$ with $\epsilon < L_\gamma(g) < L$ then $\gamma$ is disjoint from the $M^\epsilon_g$, the $\epsilon$-thin part of $(M,g)$.
\end{lemma}

{\bf Proof.} Given our upper and lower bounds for $L_c(g)$ and $R^c_g$, respectively,
there exists a "Margulis constant" $\epsilon_0$ depending only on $\alpha$ such that $M^{\epsilon_0}_g$ consists of tubes and rank two cusps. Furthermore for any $K>0$ we can choose an $\epsilon(K)< \epsilon_0$ such that the distance between $\del M^{\epsilon_0}_g$ and $\del M^{\epsilon(K)}_g$ is greater than $K$. The number $\epsilon(K)$ will only depend on $\alpha$ and $K$. Let $\epsilon = \epsilon(L/2)$. Therefore if $\gamma$ intersects $M^\epsilon_g$ it will be entirely contained in a component of $M^{\epsilon_0}_g$. The only closed geodesic in a component of $M^{\epsilon_0}_g$ will be a core curve of one of the tubes hence $L_\gamma(g) \leq \epsilon$. \qed{outofthin}

\begin{lemma}
\label{longderbound}
Let $\gamma$ be a closed, non-singular geodesic in $M_t$ such that $\|\omega_t(p)\| \leq K$ for all $p \in \gamma$. Then
$$| \cL'_\gamma(t)| \leq \sqrt{\frac{2}{3}} K L_\gamma (t).$$
\end{lemma}

{\bf Proof.}  Let $M^\gamma_t$ be the cover of $M_t$ associated to $\gamma$. For small values of $R$, $\gamma$ will have an embedded tubular neighborhood $U(R)$ of radius $R$ in $M^\gamma_t$. The $E$-valued 1-form $\omega_t$ lifts to $M^\gamma_t$ and let $K(R)$ be an upper bound for $\|\omega_t(p)\|$ for all $p \in U(R)$. Then $K(R)$ is continuous and $K(0) = K$. As we noted in the proof of Proposition \ref{lengthbound}, Hodgson and Kerckhoff show that
$$\int_{U(R)} \|\omega_t\|^2 \geq \frac{\sinh R} {\cosh R} \left( 2 + \frac{1}{\cosh^2 R} \right) \left(\frac{|\cL'_\gamma(t)|}{2L_\gamma(t)} \right)^2 \area(\del U(R)).$$
We also know that
$$\int_{U(R)} \|\omega_t\|^2 \leq K(R)^2 \int_{U(R)} dV =K(R)^2 \pi L_\gamma(t) \sinh^2 R$$
and
$$\area(\del U(R)) = 2\pi L_\gamma(t) \sinh R \cosh R.$$
Together this implies that
$$K(R)^2 \geq 2\left(2 + \frac{1}{\cosh^2 R}\right) \left(\frac{|\cL'_\gamma(t)|}{2L_\gamma(t)}\right)^2.$$
Taking the limit of both sides as $R \rightarrow 0$ and rearranging terms gives the desired inequality. \qed{longderbound}

\medskip

\noindent
{\bf Theorem \ref{controllengths}} {\em
For each $L>0$ there exists an $\epsilon >0$ and an $A> 0$ such that if $\gamma$ is a simple closed curve in $M$ with $L_\gamma(\alpha) \leq L$ and $L_c(\alpha) \leq \epsilon$ for all $c \in \cC$ then
$$e^{-AL_{\cC}(\alpha)} L_\gamma(\alpha) \leq L_\gamma(t) \leq e^{AL_{\cC}(\alpha)} L_\gamma(\alpha)$$
and
$$(1-AL_{\cC} (\alpha))\Theta_\gamma(\alpha) \leq \Theta_\gamma(t) \leq (1+AL_{\cC} (\alpha))\Theta_\gamma(\alpha)$$
for all $t \leq \alpha$.}

\medskip

\medskip

{\bf Proof.}  We begin by noting that if $L_\gamma(t)\leq e^{-4\alpha \ell_0}\ell_0$ for any $t\leq\alpha$ then the theorem follows from Proposition \ref{lengthbound} with $\epsilon = \ell_0$ and $A=4L_\cC(\alpha)$. Therefore we will assume for the remainder of the proof that $L_\gamma(\alpha) \leq L$ but $L_\gamma(t) > e^{-4\alpha \ell_0}\ell_0$ for all $t \leq \alpha$.

By Lemma \ref{outofthin} there exists a $\delta > 0$ such that if $\delta < L_\gamma(t) \leq 2L$ then the geodesic representative $\gamma_t$ of $\gamma$ in $M_t$ is disjoint from $M^{\delta}_t$. We assume that $\delta \leq e^{-4\alpha \ell_0}\ell_0$.

We need to bound the pointwise norm of $\omega_t$ on the geodesic $\gamma_t$. To do so we bound the $L^2$-norm of $\omega_t$ using Proposition \ref{normbound} and then apply the mean value theorem developed in the appendix. 

Choose $\epsilon_1>0$ such that if $L_c(t)<\epsilon_1$ then $R^c_t>3L$. Then if $L_\gamma(t)<2L$, $\gamma_t$ will not intersect the radius $R^c_t/3$ tube about $c$. In fact any ball $B_\delta$ centered at a point $p$ on $\gamma_t$ will not intersect this tube. Furthermore, since $R^c_t > \sinh^{-1} \sqrt{2}$ the area of the boundary of this tube will be universally bounded below and therefore by Proposition \ref{normbound} there exists a $K_1$ such that
$$\int_{B_\delta}\|\omega_t\|^2\leq(K_1L_\cC(\alpha))^2.$$

Let 
$$K = \frac{3 \sqrt{2 \vol(B_\delta)}K_1}{2\pi (\cosh (\delta) \sin (\sqrt{2} \delta) -\sqrt{2} \sinh (\delta) \cos (\sqrt{2} \delta))}.$$
Then the norm bound and Theorem \ref{harmboundstrain} imply that
$$\|\omega_t(p)\| \leq K L_{\cC}(t)$$
and therefore
$$|\cL'_\gamma (t)| \leq \sqrt{\frac{2}{3}} KL_{\cC}(\alpha) L_\gamma (t)$$
if $\delta < L_\gamma(t) \leq 2L$.
Next we choose $\epsilon_2$ such that 
$$e^{\sqrt{\frac{2}{3}}K \epsilon_2} \leq 2$$
and let $\epsilon = \min \{\epsilon_1, \epsilon_2\}$ and $A = \sqrt{\frac{2}{3}} K$. The rest of the argument is a repeat of the proof of Proposition \ref{lengthbound}. In particular we can first show that $L_\gamma(t) \leq 2L$ for all $t \leq \alpha$ and then integrate the derivative bound to get the final estimate.
\qed{controllengths}

\section{Bounding the Schwarzian derivative}
\label{bounding}
Recall that $\Sigma^i_t$ is the projective boundary of $M_t$ corresponding to a component of $\del N$.
The Hodge forms $\omega_t$ extend to a holomorphic quadratic differential $\Phi^i_t$ on $\Sigma^i_t$. We now use our bound on the $L^2$-norm of $\omega_t$ to bound the $L^\infty$-norm of $\Phi^i_t$. To do so we first need a local result: given a bound on the $L^2$-norm of a Hodge form on a half space $H$ we bound the $L^\infty$-norm of the quadratic differential on the projective boundary, a round disk $D$.

We  begin with a lemma from complex analysis.

\begin{lemma}
\label{schwarzbound}
Let $\Phi$ be a holomorphic quadratic differential on a round disk $D$. Then
$$\|\Phi\|_2 \geq 2\sqrt{\frac{\pi}{3}}\|\Phi(z)\|$$
for all $z \in D$.
\end{lemma}

{\bf Proof.} Without loss of generality we assume that $D$ is the unit
disk $\Delta$ in $\cx$ and $z=0$. The hyperbolic metric on $\Delta$ is $\sigma(z) =
\frac{2}{1 - |z|^2}$ and the area form in polar coordinates is $\sigma^2 rdr d\theta$. Then $\Phi = \phi dz^2$ where $\phi$ is a holomorphic function on $\Delta$. Let $\phi(z)  = \underset{n=0}{\overset{\infty}{\Sigma}} a_n z^n$ be the Taylor series for $\phi$ so $\|\Phi(0)\| = \frac{|a_0|}{4}$. We then calculate
\begin{eqnarray*}
\|\Phi\|^2_2 & = & \int_\Delta |\phi|^2 \sigma^{-4} dA \\
&= & \underset{n,m}{\Sigma} \int_0^1 \int_0^{2\pi} a_n \overline{a_m} r^{n+m} e^{\imath (n-m) \theta} \sigma^{-2} rd\theta dr \\
& = & \frac{\pi}{2}\underset{n}{\Sigma} |a_n|^2 \int_0^1 (1-r^2)^2 r^{2n+1} dr \\
& \geq & \frac{\pi}{2} |a_0|^2\int_0^1 (1-r^2)^2 rdr \\
& = & \frac{\pi}{12} |a_0|^2 \\
& = & \frac{4\pi}{3} \|\Phi(0)\|^2.
\end{eqnarray*}
\qed{schwarzbound}

Next we construct an extension of a holomorphic quadratic differential $\Phi$ on $D$ to an $E$-valued 1-form $\omega_\Phi$ on $H$. We will use a non-standard model for $\hthree$.
Namely, let $\hthree = D \times \reals$ with the Riemannian metric $\sigma^2 \cosh^2 tdx^2 + \sigma^2 \cosh^2 t dy^2 + dt^2$ where $\sigma$ is the hyperbolic metric on $D$. Then $H = D \times [0, \infty)$ is then a half space in $\hthree$ whose projective boundary is naturally identified with $D$.  We also fix an orthonormal framing of $\hthree$ by letting $\omega^1 = \sigma \cosh t dx$, $\omega^2 = \sigma \cosh t dy$ and $\omega^3 = dt$. 
Let $e_1,e_2$ and $e_3$ be the vector fields dual to the
$\omega^i$. Define $E_i$ to be the lift of $e_i$ to $\Re E$.

The holomorphic quadratic differential $\Phi$ is the Schwarzian derivative of a conformal vector field $v =f\ddz$ where $f$ is a holomorphic function with $\Phi = f_{zzz} dz^2$. To extend $\Phi$ we first use $v$ to construct a section of $E$. At each point $w = x + \imath y$ in $D$ let
$$s_\infty(w) = \left[f(w) + f_z(w)(z-w) + \frac{f_{zz}(w)}{2}(z-w)^2
\right] \ddz$$
be the projective vector field that best approximates $v$ at $w$. A projective vector field extends to a
Killing field in $\hthree$ so the function $s(w,t) = s_\infty(w)$ is an $E$-valued section on $H$. The vector fields $\Re s$ and $-\Im s$ on $H$ will extend continuously to $v$ and $\imath v$ on $D$ so $\omega_\Phi = ds$ extends continuously to $\Phi$. We will show that $\omega_\Phi$  minimizes the $L^2$-norm among all Hodge forms that extend continuously to $\Phi$.

Let ${\sf r}_w = \{w\} \times \reals$ be a geodesic in $\hthree$. We need to evaluate the extension of the parabolic vector field $\frac{(z-w)^2}{2} \ddz$ to a Killing field on the geodesic ${\sf r}_w$.

\begin{lemma}
\label{parabolicext}
The $E$-valued section $-\sigma^{-1} e^{-t}(E_1 + \imath E_2)$ evaluated on the geodesic ${\sf r}_w$ is a Killing field which extends to the projective vector field $\frac{(z-w)^2}{2} \ddz$ on $D$.
\end{lemma}

{\bf Proof.} Let $v$ be the Killing field in $\hthree$ that continuously extends to $\frac{(z-w)^2}{2}$ on $D$. Then $v$ restricted to ${\sf r}_w$ will be $-\sigma^{-1} e^{-t} e_1$. We also know that $\curl v$ is the Killing field that continuously extends to the projective vector field $\frac{i(z-w)^2}{2} \ddz$ on  $D$. From this we see that $\curl v$ restricted to ${\sf r}_w$ is the vector field $\sigma^{-1} e^{-t} e_2$.

The section $-\sigma^{-1} e^{-t}(E_1 + \imath E_2)$ evaluated at any point $p$ is the unique Killing field $w$ with $w(p) = -\sigma^{-1} e^{-t} e_1$ and $(\curl w)(p) = \sigma^{-1} e^{-t} e_2$. Therefore for $p \in {\sf r}_w$, $w = v$ proving the lemma. \qed{parabolicext}

Next we calculate the pointwise and $L^2$ norms of $\omega_\Phi$

\begin{lemma}
\label{normcalc}
If $\Phi$ is a holomorphic quadratic differential on $D$ then $\omega_\Phi$ is a Hodge form.
Furthermore the pointwise norm of $\omega_\Phi$ is
\begin{equation}
\label{pointwiseomeganorm}
\|\omega_\Phi(w,t)\| = 2 e^{-t} \sech t \|\Phi(w)\|
\end{equation}
while the $L^2$ norm on $H$ is
\begin{equation}
\label{omeganorm}
\int_H \| \omega_\Phi \|^2 =2 \|\Phi\|^2_2.
\end{equation}
\end{lemma}

{\bf Proof.} On $D$, $\Phi = \phi dz^2$ where $\phi$ is a holomorphic function. Since $s(w,t) = s_\infty(w)$,
$$\omega_\Phi(w,t) = ds(w,t) = ds_\infty(w) = \left[ \phi(w)\frac{(z-w)^2}{2} \ddz \right] dw.$$
We need to rewrite this expression in terms of $E_i$ and $\omega^i$. We first note that $dw = dx + \imath dy = \sigma^{-1}\sech t (\omega^1 + \imath \omega^2).$ From Lemma \ref{parabolicext} we know that
$$\frac{(z-w)^2}{2} \ddz = -\sigma^{-1} e^{-t}(E_1 + \imath E_2).$$
Together this implies that
\begin{equation}
\label{omegacalc}
\omega_\Phi(w,t) = -\sigma^{-2}e^{-t} \sech t \phi(w)(E_1 + \imath E_2)(\omega^1 + \imath \omega^2).
\end{equation}

Next we show that $\omega_\Phi$ is a Hodge representative. First we know that $\omega_\Phi$ is closed since $d \omega_\Phi = d(ds_\Phi) = 0$. To see that $\omega_\Phi$ is co-closed we need to calculate $\delta\omega_\Phi$. This can be done using the formula for $\delta$ given in \S2. Finally, to see that $s$ is the canonical lift of a divergence free vector field we note that $\omega_\Phi= ds$ is symmetric and traceless. The result then follows from the work in \S2 of \cite{Hodgson:Kerckhoff:cone}.

To calculate the pointwise norm of $\omega_\Phi$ on $H$ we note that
$$\|(E_1 + \imath E_2)(\omega^1 + \imath \omega^2)\|^2 = 4\omega^1 \wedge \omega^2 \wedge \omega^3$$
from which it follows that
\begin{eqnarray*}
\|\omega_\Phi(w,t)\|^2 &=& 4 \sigma^{-4} e^{-2t} \sech^2 t|\phi(w)|^2 \omega^1 \wedge \omega^2 \wedge \omega^3 \\
& = & 4e^{-2t} \sech^2 t \|\Phi(w)\|^2 \omega^1 \wedge \omega^2 \wedge \omega^3.
\end{eqnarray*} 
Next we calculate the $L^2$-norm:
\begin{eqnarray*}
\int_H \|\omega_\Phi(w,t)\|^2 &=& \int_H 4e^{-2t} \sech^2 t \|\Phi(w)\|^2 \omega^1 \wedge \omega^2 \wedge \omega^3 \\
& = & \left(\int_0^\infty 4e^{-2t} dt \right) \left(\int_D \|\Phi(w)\|^2 \sigma^2 dx dy \right) \\
&=& 2\|\Phi\|^2_2.
\end{eqnarray*}
\qed{normcalc}

We can now show that $\omega_\Phi$ minimizes the $L^2$-norm.

\begin{theorem}
\label{boundonH}
Let $\omega$ be a Hodge form on $H$ that extends continuously to a holomorphic quadratic differential $\Phi$ on $D$ and assume that $\Phi$ extends to a neighborhood of $D$ in $\chat$. Then
$$\int_H \|\omega\|^2 \geq \int_H \|\omega_\Phi\|^2 \geq \frac{8\pi}{3} \|\Phi(z)\|^2$$
for all $z \in D$.
\end{theorem}

{\bf Proof.} Let $p$ be a point on $P$ and $H_t$ the set of points in $H$ that are within a distance $t$ from $p$. Then the boundary of $H_t$ consists of a disk $P_t$ in the plane $P$ and a hemisphere $S_t$. By Theorem \ref{normform}
$$2\int_{H_t} \|\omega\|^2 = \int_{P_t \cup S_t} \imath \omega \wedge \omega^\sharp.$$

We first examine the integral over $S_t$. By definition of the Hodge star operator $\imath \omega \wedge \omega^\sharp = *_\del (\imath \omega \wedge \omega^\sharp) dA$ where $*_\del$ is the star operator for $S_t$ and $dA$ the area form, while $\omega \wedge *\omega^\sharp = *(\omega \wedge *\omega^\sharp) dV$ where $dV$ is the volume form. Both $*_\del(\imath \omega \wedge \omega^\sharp)$ and $*(\omega \wedge *\omega^\sharp)$ are real valued functions with $|*_\del(\imath \omega \wedge \omega^\sharp)| \leq *(\omega \wedge *\omega^\sharp)$.
Since $\omega$ is in $L^2$ on $H$ and $dV = dAdt$ we have
$$\int_0^\infty \int_{S_t} *(\omega \wedge *\omega^\sharp) dAdt = \int_H \omega \wedge *\omega^\sharp \leq \infty$$
and therefore
$$\underset{t \rightarrow \infty}{\lim} \int_{S_t} *(\omega \wedge *\omega^\sharp)dA = 0.$$
Since $|*_\del(\imath \omega \wedge \omega^\sharp)| \leq *(\omega \wedge *\omega^\sharp)$ we have
$$\lim_{t \rightarrow \infty} \int_{S_t} |*_\del(\imath \omega \wedge \omega^\sharp)|dA = 0$$
and
$$\int_H \|\omega\|^2 = \lim_{t \rightarrow \infty} \int_{P_t} \imath \omega \wedge \omega^\sharp.$$

We now calculate the integral over the disk $P_t$. To do so we decompose the Hodge form $\omega$ as the sum of the model Hodge form $\omega_\Phi$ and a correction term $d\tau$ where $\tau$ is the canonical lift of a vector field $w$ such that both $w$ and $\curl w$ extend to the zero vector field on $D$.
Therefore
\begin{eqnarray*}
\int_{P_t} \imath \omega \wedge \omega^\sharp &=& \int_{P_t} \imath(\omega_\Phi + d\tau) \wedge (\omega_\Phi + d\tau)^\sharp\\
 & = & 
\int_{P_t} \imath \omega_\Phi \wedge \omega_\Phi^\sharp \\
& & +2\int_{P_t} \imath d\tau \wedge \omega_\Phi^\sharp \\
& & +\int_{P_t} \imath d\tau \wedge d\tau^\sharp.
\end{eqnarray*}
By our previous remarks we know that 
$$\lim_{t \rightarrow \infty} \int_{P_t} \imath \omega_\Phi \wedge \omega_\Phi^\sharp = \int_H \|\omega_\Phi\|^2$$
and
$$\lim_{t \rightarrow \infty} \int_{P_t} \imath d\tau \wedge d\tau^\sharp = \int_H \|d\tau\|^2$$
so we need to show that
$$\lim_{t \rightarrow \infty} \int_{P_t} \imath d\tau \wedge \omega^\sharp = 0.$$

If $\alpha$ and $\beta$ are $E$-valued forms with $\alpha$ a $k$-form then
$${d}(\alpha \wedge \beta^\sharp) = d\alpha \wedge \beta^\sharp + (-1)^k \alpha \wedge (\del \beta)^\sharp$$ (On the left hand side of this equation $\alpha \wedge \beta^\sharp$ is a real valued form and $d$ is the covariant derivative for real forms.)
so
$$\int_{P_t} \imath d\tau \wedge \omega_\Phi^\sharp = -\int_{P_t} \imath \tau \wedge (\del \omega_\Phi)^\sharp + \int_{\del P_t} \imath \tau \wedge \omega_\Phi^\sharp.$$
To calculate $\del \omega_\Phi$ we recall that $d \omega_\Phi = 0$ and $\del - d = -2T$ (see \S\ref{Hyperbolic structures}) so
\begin{eqnarray*}
\del \omega_\Phi & = & -2T \omega_\Phi\\
&=& -\sigma^{-2} e^{-t} \sech t \phi (TE_1 + \imath TE_2) \wedge (\omega^1 + \omega^2) \\
&=& -\sigma^{-2} e^{-t} \sech t \phi ((E_1 +\imath E_2) \omega^1 \wedge \omega^3 + (\imath E_1 - E_2) \omega^2 \wedge \omega^3).
\end{eqnarray*}
Since there are no $\omega^1 \wedge \omega^2$ terms in $\del \omega_\Phi$,
$$\int_{P_t} \imath \tau \wedge (\del \omega_\Phi)^\sharp = 0.$$

To finish the proof we need to calculate the boundary term. To do so we use the Poincare disk model for $P$ with $p$ the center of the disk. In particular we reset our coordinates identifying the unit disk $\Delta$ in $\cx$ with the hyperbolic plane $P$. The conformal factor of the hyperbolic metric on $\Delta = P$ is $\sigma(z) = \frac{1}{1-|z|^2}$. In these coordinates (\ref{pointwiseomeganorm}) becomes
$$\|\omega_\Phi(z)\|^2 = 4\|\Phi(z)\|^2 = 4\sigma(z)^{-4}|\phi(z)|^2$$
where $\phi$ is a holomorphic function on a neighborhood of $\Delta$. In particular, $|\phi|$ is bounded on $\Delta$.

We also know that $\tau$ is the canonical lift of a vector field $w$ such that $w$ and $\curl w$ limit to zero on $\del P$. Therefore the Euclidean length of these vector fields must decay to zero at $\del P$; i.e $\sigma^{-1}\|w\|$ and $\sigma^{-1}\|\curl w\|$ extend to the zero function on $\del P$, where $\|w\|$ and $\|\curl w\|$ are hyperbolic lengths. So although $\|\tau\|^2=*(\tau \wedge *\tau^\sharp) = \|w\|^2 + \|\curl w\|^2$ may not decay to zero, there is a function $g$ on $P$ that does extend continuously to zero on $\del P$ with $\|\tau\| = \sigma g$.

Putting this all together we have
\begin{eqnarray*}
\int_{\del P_t} \imath \tau \wedge \omega_\Phi^\sharp &\leq& \int_{\del P_t} \|\tau\| \|\omega_\Phi\| dR \\
& \leq & \int_{\del P_t} 2\sigma^{-1} g|\phi| dR \\
& \leq & \int_0^{2\pi} 2g|\phi| d \theta
\end{eqnarray*}
and since $g|\phi|$ limits to zero on $\del P$ we have
$$\underset{t \rightarrow \infty}{\lim} \int_{\del P_t} \imath \tau \wedge \omega_\Phi^\sharp  = 0.$$

Therefore
$$\int_H \|\omega\|^2 = \int_H \|\omega_\Phi\|^2 + \int_H \|d\tau\|^2 \geq 2\|\Phi\|^2_2 \geq \frac{8\pi}{3} \|\Phi(z)\|$$
for all $z \in D$.
\qed{boundonH}

With this local result in place we can now bound the norm of $\Phi^i_t$. The projective structures $\Sigma^i_t$ have a fixed conformal structure $X^i$. Let $\kappa$ be the injectivity radius for the hyperbolic metric on $X^i$.

\begin{theorem}
\label{infbound}
Assume that the tube radius of $\cC$ in $M_\alpha$ is greater than $\sinh^{-1} \sqrt{2}$ and that $L_{\cC}(\alpha) \leq \ell_0$.
Then
\begin{equation}
\label{infprojbound}
L_{\cC}(t) \geq 2\sqrt{\frac{2\pi}{3}} \frac{\tanh^2(\kappa/2)}{1 +  2\|\Sigma_t\|_\infty} \|\Phi^i_t\|_\infty.
\end{equation}
\end{theorem}

{\bf Proof.} By Proposition \ref{embeddisks}, for each $z \in \Sigma_t$
there exists an embedded round disk $D$ with hyperbolic metric $\sigma_D$ such that
$$\sigma_D(z) \leq \sigma(z) \coth(\kappa/2) \sqrt{1 +
  2\|\Sigma_t\|_\infty}.$$
Therefore if we compare the norm of a holomorphic quadratic differential $\Phi$ in the $\sigma_D$-metric to the norm in the $\sigma$-metric we get
$$\|\Phi(z)\|_{\sigma_D} \geq \frac{\tanh^2(\kappa/2)}{1 + 2\|\Sigma_t\|_\infty} \|\Phi(z)\|_\sigma.$$
Let $H$ be the half space bounded by $D$ on $\Sigma_t$. By Theorem
\ref{tubebound}, $H$ is disjoint from $U^{\cC}_t$ and therefore
$$L_{\cC}(t)^2 \geq \int_{M'_t} \|\omega_t\|^2 \geq \int_{H}
\|\omega_t\|^2.$$
By Theorem \ref{boundonH}
$$\int_H \|\omega_t\|^2 \geq \frac{8\pi}{3}\|\Phi_t(z)\|_{\sigma_D}^2.$$
Combing these three formulas we have:
\begin{eqnarray*}
L_{\cC}(t)^2 & \geq & \int_H \|\omega_t\|^2 \\
& \geq & \frac{8\pi}{3} \|\Phi_t(z)\|_{\sigma_D}^2 \\
& \geq & \frac{8\pi}{3} \frac{\tanh^4(\kappa/2)}{(1+2\|\Sigma_t\|_\infty)^2} \|\Phi_t(z)\|_\sigma^2
\end{eqnarray*}
which implies the result. \qed{infbound}

Next we use the bound on $\|\Phi^i_t\|_\infty$ to bound the distance between projective structures $\Sigma^i_\alpha$ and $\Sigma^i_t$.
\medskip

\noindent
{\bf Theorem \ref{bigschwarzbound}}
{\em There exists a $C$ depending only on $\alpha$, the injectivity radius of the unique hyperbolic metric on $X^i$ and $\|\Sigma^i_\alpha\|_\infty$ such that
$$d(\Sigma^i_\alpha, \Sigma^i_t) \leq CL_{\cC}(\alpha)$$
for all $t \leq \alpha$.}

\medskip

{\bf Proof.} We will integrate (\ref{infprojbound}). Let $\sigma(T)$ be the length, in $P(X)$, of the path of projective structures $\Sigma_t$ with $t \in [T, \alpha]$. Since $\Sigma_t$ is a smooth path in $P(X)$, $\sigma(t)$ will be a smooth function and by definition $\|\Sigma_t\|_\infty \leq \|\Sigma_\alpha\|_\infty + \sigma(t)$. By Proposition \ref{schwarzpathbound}, $-\frac{d\sigma}{dt}(t) = \|\Phi_t\|_\infty$. By Proposition \ref{conelength}, $L_{\cC}(\alpha) \geq L_{\cC}(t)$ for all $t \leq \alpha$. Then by (\ref{infprojbound})
$$-\frac{d\sigma}{1/2 + \|\Sigma_\alpha\|_\infty + \sigma} \leq KL_{\cC}(\alpha)dt$$
where $K = 2\sqrt {\frac{2\pi}{3}} \coth^2 (\kappa/2)$.
Integrating both sides we have:
\begin{eqnarray*}
\int_T^\alpha -\frac{d\sigma}{1/2 + \|\Sigma_\alpha\|_\infty + \sigma} &\leq &\int_T^\alpha KL_{\cC}(\alpha) dt \\
\log \left(\frac{1/2 + \|\Sigma_\alpha\|_\infty + \sigma(T)}{1/2 + \|\Sigma_\alpha\|_\infty + \sigma(\alpha)} \right) &\leq& (\alpha - T)KL_{\cC}(\alpha) \leq \alpha KL_{\cC}(\alpha) \\
1 + \frac{\sigma(T)}{1/2 + \|\Sigma_\alpha\|_\infty} &\leq& e^{\alpha KL_{\cC}(\alpha)} \\
\sigma(T) &\leq& (1/2 + \|\Sigma_\alpha\|_\infty)\left(e^{\alpha K L_{\cC}(\alpha)} - 1\right).
\end{eqnarray*}
There exists a $C'$ depending only on $\alpha$ and $K$ (and hence $\kappa)$ such that $e^{\alpha K L_{\cC}(\alpha)} - 1 \leq C' L_{\cC}(\alpha)$. Therefore
$$\sigma(t) \leq C' L_{\cC}(\alpha)(1/2 + \|\Sigma_\alpha\|_\infty).$$
Since $d(\Sigma_\alpha, \Sigma_t) \leq \sigma(t)$ we have:
$$d(\Sigma_\alpha, \Sigma_t) \leq C L_{\cC}(\alpha)$$
where $C  = C'(1/2 + \|\Sigma_\alpha\|_\infty)$ depends only on $\kappa$, $\alpha$ and $\|\Sigma_\alpha\|_\infty$.
\qed{bigschwarzbound}

\begin{cor}
\label{projectiveconverges}
The projective structures $\Sigma^i_t$ converge to a projective structure $\Sigma ^i_{\alpha'}$ as $t \rightarrow \alpha'$.
\end{cor}

\section{Geometric limits}
We know that the projective boundary of $M_t$ converges at $\alpha'$. Now we need to show that the entire cone-manifold converges. We will need to examine geometric limits of hyperbolic cone-manifolds. Our approach will follow that of \cite{Hodgson:Kerckhoff:dehn}.

If $X$ and $Y$ are metric space and $f:X \longrightarrow Y$ is a map define $\lip(f)$ to be the infimum of all $K$ such that $d_Y(f(x_1), f(x_2)) \leq K d_X(x_1, x_2)$. The {\em bi-Lipschitz distance} between $X$ and $Y$ is $\bilip(X,Y) = \inf \{|\log \lip(f)| + |\log \lip(f^{-1})| | f:X \longrightarrow Y \mbox{\ is bi-Lipschitz}\}$. This defines the {\em bi-Lipschitz topology} on the set of metric spaces.

To show that a sequence of compact hyperbolic manifolds with boundary converges we need to control three quantities: the prinicipal curvatures of the boundary, the injectivity radius and the width of collar neighborhoods of the boundary. If $M$ is a hyperbolic manifold we let $\inj_M = \inf \{\inj_x | x \in M\}$. We define $\del M(t)$ to be those points in $M$ whose distance from $\del M$ is less than $t$. Then $\width(\del M) = \inf\{t | \del M(t') \mbox{\ is a collar of $M$ for all $t' < t$}\}$. The geometric convergence theorem that follows is essentialy due to Kodani \cite{Kodani:convergence}(see the remarks on p. 20 of \cite{Hodgson:Kerckhoff:dehn}):

\begin{theorem}
\label{geomboundaryconv}
Let $\lambda^-$, $\lambda^+$, $i_0$ and $W$ be real constants with $\lambda^+ \geq \lambda^-$ and $i_0, W >0$. Let $M_n$ be a sequence of hyperbolic manifolds with boundary such that the principal curvatures of $\del M_n$ are contained in the interval $[\lambda^-, \lambda^+]$, $\inj_{M_n} \geq i_0$ and $\width(\del M_n) \geq W$. Then there exists a hyperbolic 3-manifold with boundary $M_\infty$ and a subsequence $\{n_k\}$ such that $M_{n_k} \rightarrow M_\infty$ in the bi-Lipschitz topology. Furthermore if all the $M_n$ are diffeomorphic and have bounded volume (or diameter) then $M_\infty$ is diffeomorphic to $M_n$.
\end{theorem}

We will apply this result to suitably chosen compact submanifolds of our hyperbolic cone-manifolds $M_t$.

\subsection{Geometric limit of geometrically finite ends}
\label{univalent}
To construct these compact submanifolds we remove a neighborhood of the geometrically finite ends. To do this we need to understand how the projective boundary determines the
hyperbolic geometry of the geometrically finite ends. This information comes from work of Epstein (\cite{Epstein:horospheres}) and Anderson (\cite{Anderson:projective}) which we will review here. 

We will use the same coordinates for $\hthree$ as we did in \S \ref{bounding}. In particular, let $U$ be the upper half space of $\cx$ with hyperbolic metric $\sigma$. Then $\hthree=U\times\reals$ with metric $$\sigma^2 \cosh^2tdx^2 + \sigma^2 \sinh^2 t dy^2 +dt^2.$$
Let $P_d$ be the set of points of the form $(z,d)$ in $\hthree$. Then $P_0$ is a hyperbolic plane and $P_d$ is a constant curvature plane a (signed) distance $d$ from $P_0$.

Let $\psi:U \longrightarrow \chat$ be a locally univalent map and $\Phi = S\psi$ its Schwarzian derivative. The {\em
osculating M\"obius} transformation $M_{\psi(z)}$ is the unique M\"obius
transformation whose two jet agrees with $\psi$ at $z$. We define $\Psi: \hthree \rightarrow \hthree$
by
$$ \Psi(z,d)= M_{\psi(z)}(z,d).$$
Note that $\Psi$ extends continuously to $\psi$ on $U$.

The following two results can be found in \S3 of \cite{Anderson:projective}.
\begin{prop}
\label{epder}
Let $p = (z,d)$ be a point in $\hthree$. There exist
an orthonormal basis $\{e_1, e_2, e_3\}$ for $T\hthree_p$ with $e_1$
and $e_2$ spanning the plane normal to $P_d$ and an orthonormal
basis for $T\hthree_{\Psi(p)}$ such that $d\Psi$ at $p$ in these coordinates
is:
$$\left( \begin{array}{ccc}
1 + \frac{\|\Phi(z)\|}{4e^d \cosh d} & 0 & 0\\
0 & 1 - \frac{\|\Phi(z)\|}{4e^d \cosh d} & 0\\
0 & 0 & 1
\end{array} \right).$$
In particular if $4e^d \cosh d > \|\Phi\|_{\infty}$ then $\Psi$ is an
orientation preserving local diffeomorphism at $p$.
\end{prop}

{\bf Remark.} For each $d$ the nearest point retraction $\pi_d$ defines a natural map from $U$ to $P_d$. In \cite{Anderson:projective} what is actually calculated is the derivative of the composition $\Psi \circ \pi_d$. Since the derivative of $\pi_d$ is easy to calculate we can translate the work in \cite{Anderson:projective} to the above proposition.\footnote{There is an error in the calculation of the eigenvalues on p. 35 of \cite{Anderson:projective}. They should be $\frac{1}{2}(1+1/t) + \|Sf(0)\|/4t$ and $\frac{1}{2}(1+1/t) - \|Sf(0)\|/4t$ not $\frac{1}{2}(1+1/t) + \|Sf(0)\|$ and $\frac{1}{2}(1+1/t) - \|Sf(0)\|$.}

When $\Psi$ is an immersion there are also formulas for the curvature of the image surface $\Psi(P_d)$.

\begin{prop}
\label{princcurv}
Let $p = (z,d)$ be a point on $P_d$ and let $k_1 =
-\frac{\|\Phi(z)\|}{\|\Phi(z)\| - 1}$ and $k_2 =
-\frac{\|\Phi(z)\|}{\|\Phi(z)\| + 1}$. Then the
principal curvatures (if defined) of $\Psi(P_d)$ at $\Psi(p)$ are
$$\frac{\sinh d + k_i \cosh d}{\cosh d + k_i\sinh d}$$
for $i=1,2$.
\end{prop}

Let $\Sigma$ be a projective structure on a surface $S$ with conformal structure $X$ and let $\Sigma_F$ be the fuchsian projective structure with conformal structure $X$. Then there is a representation $\rho_F: \pi_1(S) \longrightarrow PSL_2\reals$ such $\Sigma_F = U/\rho_F(\pi_1(S))$. Identifying $U$ with the universal cover $\tilde{S}$ and the deck transformations with $\rho_F(\pi_1(S))$ there is a conformal developing map $\psi:U \longrightarrow \chat$ for $\Sigma$. In particular, $\Sigma$ has a holonomy representation $\rho$ and $\psi \circ \rho_F(\gamma) = \rho(\gamma) \circ \psi$ for all $\gamma \in \pi_1(S)$.
As above $\psi$ extends to a map $\Psi: \hthree \longrightarrow \hthree$. It is clear from the definition that this construction is natural. That is $\Psi \circ \rho_F(\gamma) = \rho(\gamma) \circ \Psi$ for all $\gamma \in \pi_1(S)$.

The group $\rho_F(\pi_1(S))$ also acts on $\hthree$ with quotient homeomorphic to $S \times \reals$. 
We can therefore view $\hthree$ as the universal cover of $S \times \reals$, identifying $\tilde{S} \times \{d\}$ with $P_d$ in $\hthree$.  Then $\Psi$ is a map from $\tilde{S}  \times \reals$ to $\hthree$. Restricted to $\tilde{S} \times [d, \infty)$ where $e^d > \sqrt{2\|\Sigma\|_\infty + 1}$, $\Psi$ is a diffeomorphism. Let $\cE(\Sigma, d)$ be the hyperbolic structure on $S \times [d,\infty)$ defined by this developing map. The hyperbolic structure $\cE(\Sigma, d)$ extends to the projective structure $\Sigma$ on $S \times \{\infty\}$ so $\cE(\Sigma, d)$ is a geometrically finite end with projective boundary $\Sigma$. The plane $P_t$ descends to surfaces $\cS(\Sigma, t)$ that foliate $\cE(\Sigma, d)$.

\begin{prop}
\label{epconvex}
The surfaces $\cS(\Sigma, t)$ are convex in $\cE(\Sigma,d)$ and are strictly convex if $t > 0$.
\end{prop}

{\bf Proof.} This is a direct consequence of Proposition \ref{princcurv}. \qed{epconvex}

The foliation of $\cE(\Sigma, d)$ by convex surfaces implies that $\cE(\Sigma, d)$ embeds in a hyperbolic cone-manifold with projective boundary $\Sigma$. More precisely:

\begin{prop}
\label{epembed}
If $\Sigma$ is a component of the projective boundary of a
hyperbolic cone-manifold $M$ then $\cE(\Sigma, d)$ embeds in $M$ if $d>0$.
\end{prop}

{\bf Proof.} The proof is the same as Proposition \ref{allhalfembed}. \qed{epembed}

Although $\cE(\Sigma, d)$ embeds it may intersect the tubes $U^c_t$. We need to show that $d$ can be chosen large enough so that this doesn't happen. To do this we will use an alternative construction of the surfaces $\cS(\Sigma, d)$ as an envelope of horospheres.

For each $p \in \hthree$ the identification of the unit sphere in $T_p \hthree$ with the ideal boundary $\chat$ of $\hthree$ determines a {\em visual measure} $\mu_p$ of $\chat$. Given a conformal metric $\sigma$ on a domain $\Omega \subseteq \chat$ and a point $z \in \chat$ the set of points in $\hthree$ whose visual measure equals $\sigma$ at $z$ is a horosphere $\cH_{z,\sigma}$. This horosphere also has the following property: A plane $P$ in $\hthree$ limits to a round circle in $\chat$. This circle bounds two disks and we assume one of these disks $D$ contains $z$. Then hyperbolic metric $\sigma_D$ on $D$ will agree with $\sigma$ at $z$ if and only if $P$ is tangent to $\cH_{z,\sigma}$.

The envelope of this family of horospheres $\cH_{z,\sigma}$ is a surface in $\hthree$. A similar construction works in a geometrically finite end. In fact if $\sigma$ is the hyperbolic metric on $\Sigma$ and $\sigma_d = e^d \sigma$ the envelope of the family of horospheres $\cH_{z,\sigma_d}$ is the surface $\cS(\Sigma, d)$.

\begin{theorem}
\label{enddisjoint}
Let $\kappa$ be the injectivity radius of the hyperbolic metric on $\Sigma$.
If $e^d > \coth(\kappa/2)\sqrt{1 + 2\|\Sigma\|_\infty}$ then $\cE(\Sigma, d)$ is embedded and disjoint from $U^{\cC}_{t}$.
\end{theorem}

{\bf Proof.} Since the envelope of the family of horospheres $\cH_{z, \sigma_d}$ is the surface $\cS(\Sigma, d)$ the union of the horoballs bounded by $\cH_{z, \sigma_d}$ is the entire end $\cE(\Sigma, d)$. By Proposition \ref{embeddisks} there exists an embedded round disk $D'$ with hyperbolic metric $\sigma_D$ such that
$$\sigma_{D'}(z) < \sigma(z) \coth(\kappa/2)\sqrt{1+2\|\Sigma\|_\infty}.$$
By our choice of $d$ there exists a round disk $D \subset D'$ such that
$\sigma_D(z) = \sigma_d(z)$. The round disk $D$ is the projective boundary of an embedded halfspace $H$ which will have hyperbolic boundary a plane $P$. The plane $P$ will be tangent to $\cH_{z, \sigma_d}$ so the horoball bounded by $\cH_{z,\sigma_d}$ will be contained in $H$. By Theorem \ref{tubebound} $H$ is disjoint from $U^{\cC}_t$ and therefore the horoball bounded by $\cH_{z, \sigma_d}$ is disjont from $U^{\cC}_t$ proving the theorem. \qed{enddisjoint}


For each $\Sigma$ the map $\Psi$ canonicaly identifies $\cE(\Sigma, d)$ as a Reimannian metric on a fixed copy of $S \times [d, \infty)$. We the have the following proposition:

\begin{prop}
\label{endsconv} Let $\Sigma_t$ be a sequence of projective structures which converge to $\Sigma_\infty$ in $P(X)$ and assume that $e^d > \sqrt{2\|\Sigma_t\|_\infty + 1}$ for all $t$. Then the metrics $\cE(\Sigma_t, d)$ converge to $\cE(\Sigma, \infty)$ in the compact-$C^\infty$ topology on metrics on $S \times [d, \infty)$.
\end{prop}

{\bf Proof.} This is a simple consequence of the fact that the maps $\tilde{\Psi}_t$ depend continuously on the projective structures $\Sigma_t$. In particular, if the $\Sigma_t$ converge to $\Sigma_\infty$ then $\tilde{\Psi}_t$ can be chosen to converge to $\tilde{\Psi}_\infty$ in the $C^\infty$-topology on maps from $\tilde{S} \times [d, \infty)$ to $\hthree$ which implies that the metrics converge.
\qed{endsconv}

On a compact manifold converge in the compact-$C^\infty$ topology of a sequence of metrics implies that the associated metric spaces converge in the bi-Lipschitz topology. Since $S \times [d, \infty)$ is non-compact we only get bi-Lipschitz convergence on compact submanifolds such as collars $S \times [d,d']$. This will be enough for our applications.

\subsection{The Schl\"{a}fli formula}
We will need to bound the volume of the complement of the geometrically finite ends in the cone manifolds $M_t$.  To do so we will use the generalized Schl\"{a}fli
formula of Rivin and Schlenker. Although their formula applies in much
greater generality we will stick to the case of a 3-manifold with
boundary $M$ and a smooth family of hyperbolic cone-metrics $g_t$ where $t$ is the cone-angle.
Let ${\rm
  I}(t)$ and ${\rm II}(t)$ be the first and second fundamental forms
for $\del M$ in the $g_t$ metric and $H(t)$ the mean curvature. Finally let $V(t)$ be
the volume of $M$ in the $g_t$ metric. The generalized Sch\"afli
formula is then:

\begin{theorem}[Rivin-Schlenker \cite{Rivin:Schlenker:schlafli}]
\label{schafli}
$$-\frac{1}{3} V'(t) = \int_{\del M} \left(H'(t) + \frac{1}{2} \langle {\rm
    I}'(t), {\rm II}(t) \rangle \right) dA  + L_{\cC}(t)$$ 
\end{theorem}

We will use the following simple corollary of this result:

\begin{cor}
\label{schaflibound}
Assume that the metric $g_t$ are defined for $a<t\leq b$ and that as
$t \rightarrow a$, $g_t$ converges in the $C^\infty$-topology on a
collar neighborhood of $\del M$. Furthermore assume that $L_{\cC}(t)$ is bounded. Then $V(t)$ is
bounded. 
\end{cor}

{\bf Proof.} Since $g_t$ converges the quantities $H'(t)$ and ${\rm
  I}'(t)$ and ${\rm II}(t)$ are all bounded and therefore the integral 
$$\int_{\del M} \left(H' + \langle {\rm I}', {\rm II} \rangle \right) dA$$
will be bounded. Since $L_{\cC}$ is bounded, Theorem \ref{schafli} implies that $V'$ is bounded which in turn
implies that 
$$V(T) = \int_1^T V' dt + V(1)$$
is bounded. \qed{schaflibound}

{\bf Remark.} This result could also have been proven using the standard Schl\"{a}fli formula for manifolds with polyhedral boundary. One simply needs to construct a polyhedral approximation for the smooth boundary.

\subsection{Geometric limits of cone-manifolds}
By Corollary \ref{projectiveconverges} we know that the projective structures $\Sigma^i_t$ converge to a projective structure $\Sigma^i_{\alpha'}$ as $t \rightarrow \alpha'$. Therefore there exists a $d>0$ such that $e^{d-1} > \sqrt{2\|\Sigma^i_t\|_\infty + 1}$ for all $t \in (\alpha', \alpha]$ and $i=1, \dots, n$. We then let $\cM_t$ be the closure of $M'_t \backslash \left( \overset{n}{\underset{i=1}{\bigcup}} \cE(\Sigma^i_t, d) \right)$. The boundary of $\cM_t$ consist of the boundary $\del U^c_t$ of the tubes and the boundary $\cS(\Sigma^i_t, d)$ of the geometrically finite ends.

\begin{theorem}
\label{coneconv}
There exists a sequence $\{t_n\}$ in $(\alpha', \alpha]$ with $t_n
\rightarrow \alpha'$ such that $\cM_{t_n}$ converges in the
bi-Lipschitz topology to a hyperbolic cone-manifold $\cM_{\alpha'}$ homeomorphic to $\cM_t$. 
\end{theorem}

{\bf Proof.} We need to see that the conditions of Theorem \ref{geomboundaryconv} hold for the family $\cM_t$. As we have already noted the norms $\|\Sigma^i_t\|_\infty$ are bounded since the projective structures converge. Proposition \ref{princcurv} then implies that the principal curvatures of $\cS(\Sigma^i_t, d)$ are bounded above and below. Since the radii of the tubes $U^c_t$ are greater that $\sinh^{-1}\sqrt{2}$ the principal curvatures of $\del U^c_t$ are bounded between $\coth (\sinh^{-1} \sqrt{2})$ and $\tanh (\sinh^{-1} \sqrt{2})$. We bounded the injectivity radius in Proposition \ref{injbound}. Our choice of $d$ gaurantees that $\width(\cS(\Sigma^i_t, d)) \geq 1$. To see that the boundary components $\del U^c_t$ have a definite width we note that the tubes $U^c_t$ are not maximal. In fact from Proposition \ref{firsttubebound} we see that there are disjoint embedded tubes whose torus boundary has area $0.51$. The difference between these tubes and the $U^c_t$ will be a collar of definite width for all $t$.

This is enough to obtain a limiting hyperbolic manifold $\cM_{\alpha'}$. To see that $\cM_{\alpha'}$ is homeomorphic to the $\cM_t$ we need to bound the volume of the $\cM_t$. If we consider the $\cM_t$ as a family of metrics on a fixed manifold Proposition \ref{endsconv} implies that we can choose these metrics such that they converge on a neighborhood of each $\cS(\Sigma^i_t, d)$. We then apply Corollary \ref{schaflibound} to $M_t \backslash \left( \overset{n}{\underset{i=1}{\bigcup}} \cE(\Sigma^i_t, d) \right)$ to bound the volume. \qed{coneconv}

\begin{theorem}
\label{convergetoprime}
There exists a $M_{\alpha'} \in GF(N, \cC)$ such that $M_t \rightarrow M_{\alpha'}$.
\end{theorem}

{\bf Proof.} The boundary of $\cM_{\alpha'}$ consists of tori and higher genus surfaces. On a collar of the higher genus ends the manifolds $\cM_{t_n}$ converge to a collar of the geometrically finite ends $\cE(\Sigma^i_{\alpha'}, d-1)$.
Therefore we can glue the ends $\cE(\Sigma^i_{\alpha'}, d-1)$ to the higher genus boundary components. It is shown \S3 of \cite{Hodgson:Kerckhoff:dehn} that the metric can be extended to a cone singularity (or cusp if $\alpha' = 0$) at the torus components of $\del \cM_{\alpha'}$ with cone angle $\alpha'$. In fact they show more than this. For large $n$ there are bi-Lipschitz diffeomorphisms $f_n$ from $\cM_{t_n}$ to $\cM_{\alpha'}$. For a fixed $c \in \cC$, $f_n$ maps $\del U^c_{t_n}$ to a fixed component $\del U^c_{\alpha'}$ of $\del \cM_{\alpha'}$. The {\em meridian} of $\del U^c_{t_n}$ is the unique homotopy class a of non-trivial simple closed curve on $\del U^c_{t_n}$ that bounds a disk in $U^c_{t_n}$. Hodgson and Kerckhoff further show that $f_n$ maps the meridians to a fixed homotopy class on $\del U^c_{\alpha'}$ and this homotopy class is a meridian of the cone singularity in the extended structure.

The extended manifold $M_{\alpha'}$ is then a geometrically finite cone-metric on a pair $(\hat{N}, \hat{\cC})$.  Since the maps $f_n$ take meridians to meridians the $f_n$ extend to homeomorphisms from $(N, \cC)$ to $(\hat{N}, \hat{\cC})$. These extensions of $f_n$ can also be chosen to be conformal maps form the conformal boundaries of $M_{t_n}$ to $M_{\alpha'}$. By Theorem \ref{param} in open an interval about $\alpha'$ there exists a one parameter family of cone-manifolds $\hat{M}_t$ with cone angle $t$, 
conformal boundary $X$ and $\hat{M}_{\alpha'} = M_{\alpha'}$. We need to show that each $f_n$ is homotopic to an isometry from $M_{t_n}$ to $\hat{M}_{t_n}$.

Let $\rho_t$ and $\hat{\rho}_t$ be the holonomy representations of $M_t$ and $\hat{M}_t$ respectively. (Note that they are representations of $\pi_1(N\backslash \cC) \isom \pi_1 (\hat{N}\backslash \hat{\cC})$ not of $\pi_1(N) \isom \pi_1(\hat{N})$.) Convergence in the bi-Lipschitz topology implies that the representations $(f_n)_* \rho_{t_n}$ converge to $\hat{\rho}_{\alpha'}$. By Theorem 5.7 in \cite{Bromberg:rigidity} the space of conjugacy classes of representations is locally parameterized by the complex length of the meridians and the conformal boundary. This is a stronger version of Theorem \ref{param} which allows representations where the holonomy of the meridians is not elliptic. It implies that $\hat{\rho}_{t_n} = (f_n)_* \rho_{t_n}$ for large $n$. By Theorem 1.7.1 of \cite{Canary:Epstein:Green} on $\cM_{t_n}$ $f_n$ will be homotopic to an isometric embedding of $\cM_{t_n}$ in $\hat{M}_{t_n}$. Since $\cM_{t_n}$ extends to a geometrically finite cone-manifold in a unique way this implies that $f_n$ is homotopic to an isometry from all of $M_{t_n}$ to $\hat{M}_{t_n}$.

To finish the proof we choose a fixed large value of $n$ and use the map $f_n$ to pull back metrics in $GF(\hat{N}, \hat{\cC})$ to metrics in $GF(N, \cC)$. Under this identification $M_{t_n} = \hat{M}_{t_n}$. Theorem \ref{param} then implies that $M_t = \hat{M}_t$ wherever both are defined. Therefore $M_t \rightarrow M_{\alpha'}$ as desired.

\qed{convergetoprime}

We are now ready to prove Theorem \ref{downtozero}:

\medskip

\noindent
{\bf Theorem \ref{downtozero}} {\em
Let $M_\alpha \in GF(N, \cC)$ be a geometrically finite hyperbolic
cone-metric with cone angle $\alpha$. Suppose the tube radius of the cone singularity is $\geq
\sinh^{-1} \sqrt{2}$. Then there exists an
$\ell_0$ depending only on $\alpha$ such that if $L_c(\alpha) \leq \ell_0$ for all $c \in \cC$ then the one parameter family of cone-manifolds $M_t \in GF(N, \cC)$ is defined for all $t \leq \alpha$.}

\medskip

{\bf Proof.} By Theorem \ref{param} the interval for which the family $M_t$ is defined is open in $[0, \alpha]$. By Theorem \ref{convergetoprime} this interval is also closed. Therefore $M_t$ is defined for all $t \in [0,\alpha]$.
\qed{downtozero}

\section{Rank two cusps}

In this section we show how bounds on the $L^2$-norm control the shape of a rank two cusp. Recall that a rank two cusp is the quotient of a horoball, centered at infinity in the upper half space model, by parabolic elements $z\mapsto z+1$ and $z\mapsto z+\tau$ with $\Im \tau > 0$.
The horoball is foliated by horospheres which are horizontal planes in $\hthree$. The quotient of these planes are tori which foliate the rank two cusp. Each tori will be conformally equivalent with $\tau$ the Teichm\"{u}ller parameter of the tori.
To normalize the cusp we choose the horoball so that in the quotient the boundary torus has area $\frac{1}{2}$.

We will use similar notation for cusps as we do for short geodesics and their tubular neighborhoods. In particular if $\gamma$ is a torus component of $\del N$ the $U^\gamma_t$ will be the associated rank two cusp. We also let $\cL_\gamma (t)$ be the Teichm\"{u}ller parameter of the cusp.

We do not know, a priori, that the cusps $U^\gamma_t$ are embedded. The proof of this is essentially the same as Theorem \ref{tubebound}. In particular we have:

\begin{prop}
\label{cuspsinbed}
If the tube radius of the cone-singularity is greater than $\sinh^{-1} \sqrt{2}$ and $L_c(t) \leq \ell_0$ for all $c \in \cC$ then the cusps $U^\gamma_t$ are embedded and disjoint from the tubes $U^c_t$ for all $c \in \cC$.
\end{prop}

Next we need to control the derivative of the Teichm\"{u}ller parameter  as $t$ varies. Note that $\cL_\gamma (t)$ is a point in the Teichm\"{u}ller space of a torus which is canonically identified with ${\bb H}^2$ so we will measure the derivative $\cL'_\gamma (t)$ in the hyperbolic metric. We then have

\begin{theorem}
\label{cuspderivative}
$$\cL_\cC (\alpha)^2 \geq |\cL'_\gamma (t)|^2$$
\end{theorem}

{\bf Proof.} 
The proof is similar to the proof Theorem \ref{boundonH}.

On the rank two cusp the Hodge form $\omega_t$ is a sum of a model deformation $\omega ^m_t$ and a correction term $\omega ^c_t$. The model is constricted in \S3.7 of \cite{Bromberg:rigidity}.
The term $\omega ^c_t$ is trivial so $\omega ^c_t = d\psi_t$ where $\psi_t$ is an $E$-valued section  on $U^\gamma_t$.

Recall that to calculate the $L ^2$-norm of $\omega_t$ on $U ^t _\gamma$ we need to integrate $\imath \omega_t \wedge \omega _t ^\sharp$ over the boundary torus $\del U ^t _\gamma$. As in the proof of Theorem \ref{boundonH} we expand $\imath \omega_t \wedge \omega_t ^\sharp$ to get 

$$\imath \omega _t \wedge \omega _t ^\sharp = \imath \omega ^m _t \wedge (\omega ^m _t) ^ \sharp + 2 \imath \omega ^c _t \wedge (\omega ^m _t) ^\sharp + \imath \omega ^c _t \wedge (\omega ^c _t) ^\sharp$$
Using the integration by parts argument from  Theorem \ref{boundonH} we have
$$\int _{\del U ^t _\gamma} \imath \omega ^c _t \wedge (\omega ^m _t) ^\sharp = -\int_{\del U^t _\gamma} \imath \psi _t \wedge (\del \omega^m _t) ^\sharp.$$
From the explicit form of $\omega^m_t$ we see that $\del \omega ^m _t$ has no terms tangent to $\del U^\gamma _t$ so 
$$\int_{\del U^\gamma _t} \imath \psi_t \wedge (\del \omega^m _t)^\sharp=0$$
which implies that
$$\int_{U^t_\gamma} \|\omega_t\|^2 = \int_{U_\gamma ^t} \|\omega ^m _t\| ^2 + \int _{U^t_\gamma} \|\omega^c_t\|^2.$$ 
We can also calculate the $L^2$-norm
$$\int_{U_\gamma ^t} \|\omega^m_t\|^2= |\cL'_\gamma (t)|^2.$$

By Propositions \ref{normbound} and \ref{cuspsinbed} we know that
$$\cL_\cC (\alpha)^2\geq \int_{M_t}\|\omega_t\|^2 \geq \int_{U^\gamma_t} \|\omega_t\|^2.$$
combining this last inequality with the previous two equalities completes the proof.
\qed{cuspderivative}

As an immediate corollary we have:

\begin{cor}
\label{cuspsconverge}
The length of the path $\cL_\gamma(t)$ with $t \in (\alpha',\alpha]$ is bounded and therefore $\cL_\gamma (t)$ converges to $\cL_\gamma(\alpha')$ as $t\rightarrow \alpha'$, where $\cL_\gamma(\alpha')$ is a complex number with $\Im \cL_\gamma (\alpha')>0$.
\end{cor}

The proofs of Theorems \ref{bigschwarzbound} and \ref{controllengths} do not change when rank two cusps are allowed. To prove Theorem \ref{downtozero} we need to modify the geometric limit arguement. When there are rank two cusps we define $\cM_t$ to be the complement of geometrically finite ends, the tubes about the cone-singularity and the rank two cusps. Then just as in Proposition \ref{boundinj}, Corollary \ref{cuspsconverge} implies that the injectivity radius of $\cM_t$ is bounded below. Once we have this the geometric limit arguement and the proof of Theorem \ref{downtozero} follow as before.

\section{Applications to Kleinian groups}
We now use our results to prove some estimates on pinching short curves on the conformal boundary of a smooth hyperbolic 3-manifold. We need to reset our notation. Let $M$ be a complete, smooth hyperbolic 3-manifold and $\cC$ a collection of disjoint simple closed curves on the conformal boundary $\del M$ of $M$. We also assume that the curves in $\cC$ are homotopically distinct in $M$ and that each $c \in \cC$ represents a primitive element of $\pi_1(M)$.  Then $L_c(M)$ will be the length of the geodesic representative of $c$ in $M$ and $L_c(\del M)$ will be the length of the geodesic representative of $c$ for the hyperbolic metric on the conformal boundary. We need the following preliminary result which is a combination of theorems of Canary (Theorem 5.1 in \cite{Canary:retraction}) and Otal (Theorem 3 in  \cite{Otal:unknottwo}).

\begin{theorem}
\label{unknotted}
There exists an $\epsilon_0 > 0$ such that if $L_c(\del M) \leq \epsilon_0$ for each $c \in \cC$ then the geodesic representatives $c^*$ of $c$ in $M$ is isotopic to $c$ on $\del M$. Furthermore this isotopy is disjoint from the geodesic representatives of the other curves in $\cC$.
\end{theorem}

In the next two results we compare the geometry of a hyperbolic manifold with short curves in its conformal boundary to the geometry of a hyperbolic manifold with those curves pinched to rank one cusps.

\begin{theorem}
\label{curt}
Assume $M$ is a smooth geometrically finite hyperbolic 3-manifold without rank one cusps. Assume that $\cC$ is a collection of homotopically distinct and disjoint simple closed curves on $\del M$ and that each $c \in \cC$ represents a primitive element of $\pi_1(M)$. There exists an $\ell'_0 > 0$ such that if $L_c(\del M) \leq \ell'_0$ for all $c \in \cC$ the following holds:
\begin{enumerate}
\item There exists a smooth, geometrically finite hyperbolic structure $\hat{M}$ homeomorphic to $M$ with each curve $c$ pinched to a rank one cusp.

\item The components of the conformal boundaries of $M$ and $\hat{M}$ that are disjoint from $\cC$ are isomorphic.

\item If $X$ is a component of the conformal boundary disjoint from $\cC$ and $\Sigma$ and $\hat{\Sigma}$ are the projective boundaries on $X$ for $M$ and $\hat{M}$, respectively, then
there exists a $C$ depending only on the injectivity radius of the hyperbolic metric on $X$ and $\|\Sigma\|_\infty$ such that
$$d(\Sigma, \hat{\Sigma}) \leq C L_{\cC}(\del M).$$
\end{enumerate} 
\end{theorem}

{\bf Proof.} Let $\ell'_0 = \min\{\ell_0, \epsilon_0\}$. By Theorem \ref{downtozero} there exists a one-parameter family of cone-manifolds $M_t$ for $t \in [0,2\pi]$ with cone singularity $\cC^*$ such that $M_{2\pi} = M$. By Theorem \ref{unknotted}, the manifold $M$ has a compact core $M'$ which is disjoint from $\cC^*$. Furthermore each $c^*$ will be isotopic to a curve on the boundary of $M'$. Recall that $M_0$ is complete hyperbolic structure on the topological manifold $M \backslash \cC$. Then the cover $\hat{M}$ of $M_0$ associated to the compact submanifold $M'$ is the hyperbolic manifold satisfying (1) and (2).

Finally, (3) holds by Theorem \ref{bigschwarzbound}. Note that if $M$ has incompressible boundary $C$ depends only on the injectivity radius of $X$ not on $\|\Sigma\|_\infty$. This is because when $\del M$ has incompressible boundary $\|\Sigma\|_\infty \leq 3/2$ by Nehari's Theorem.
\qed{curt}

This theorem should be compared to Corollary 1.3 in \cite{McMullen:cusps} which treats the special case where $M$ is quasifuchsian. That result has a better bound than the one we achieve here. Our bound could be improved be finding a lower bound on the distance between the tubes $U^c_t$ and an embedded half space $H$ rather than just showing that $U^c_t$ and $H$ are disjoint.

In our next we result we control the complex length of geodesics in $M$. In this theorem $\Theta_\gamma(M)$ is the imaginary part of the complex length.

\begin{theorem}
\label{saar}
For each $L>0$ there exists an $\epsilon >0$ and an $A >0$ such that if $\gamma$ is a closed geodesic in $M$ with $L_\gamma(M) \leq L$ and $L_c(\del M) \leq \epsilon$ for all $c \in \cC$ then
$$e^{-AL_{\cC}(\del M)}L_\gamma(\hat{M}) \leq L_\gamma(M) \leq e^{AL_{\cC}(\del M)}L_\gamma(\hat{M})$$
and
$$(1-AL_{\cC}(M))\Theta_\gamma(\hat{M}) \leq \Theta_\gamma(M) \leq (1+AL_{\cC}(M))\Theta_\gamma(\hat{M}).$$
\end{theorem}

{\bf Proof.} This follows immediately using the construction of $\hat{M}$ in the previous result and Theorem \ref{controllengths}. \qed{saar}

This result should be compared with Proposition 5.1 in \cite{Canary:Culler:Hersonsky:Shalen} which is essentially the same result.

\section{Appendix - Mean value inequalities}
In this appendix we prove mean value inequalities for harmonic vector and strain fields. For strain fields this inequality was proved by Hodgson and Kerckhoff in an early version of \cite{Hodgson:Kerckhoff:dehn}. Their proof is not in the current version of \cite{Hodgson:Kerckhoff:dehn} as they have found a simpler proof of their main results which does not require the inequality. With their permission we recount the result here.

We also prove a mean value inequality for vector fields $v$ where $\|\Delta v\|$ is bounded. The proof is essentially the same as for strain fields if not simpler. Will start with the vector field inequality along with an application to Hodge forms on geometrically finite ends.

Using the identification of the tangent bundle with the real part of $E$, a
vector valued $k$-form $\omega$ can be identified as a real $E$-valued $k$-form. Then $\Delta \omega$ will be the Laplacian associated to the bundle $E$. For a function $u$, $\Delta u$ will be the standard Laplacian on functions. It will be clear from context which Laplacian we are using.

\begin{lemma}
\label{boundonLu}
If $\|\Delta v\| \leq b$ then $-(\Delta \|v\| + 2\|v\|) \geq -b$.
\end{lemma}

{\bf Proof.} For vector fields we have the Weitzenb\"ock formula,
$$\Delta v  = \nabla^* \nabla v + 2v$$
where $\nabla$ is the Riemannian connection and $\nabla^*$ its adjoint (see \S2 of \cite{Hodgson:Kerckhoff:cone}).

Let
$$\nabla^2_{XY} = \nabla_X \nabla_Y - \nabla_{\nabla_X Y}.$$
Then $\nabla^* \nabla = -\Sigma_i \nabla^2_{e_i e_i}$, where $\{e_1, e_2, e_3\}$ is an orthonormal frame field on $M$. For functions $\Delta = \nabla^*\nabla$. 
Using this formula we see that for any tensor $S$ on a Riemannian manifold:
$$\Delta \|S\|^2 = 2 \langle \nabla^* \nabla S, S \rangle - 2\|\nabla S\|^2.$$
Combining this formula with the Weitzenb\"{o}ck formula for $v$ we have
\begin{equation}
\label{form1}
\Delta \|v\|^2 = 2\langle \Delta v, v \rangle - 4\|v\|^2 - 2\|\nabla v\|^2.
\end{equation}

Let $u = \|v\|^2$. Applying the product formula for the Laplacian to
$u$ times itself gives 
$$\Delta(u^2) = 2u \Delta u - 2\|\nabla u\|^2.$$
Combining this formula with (\ref{form1}) we get:
$$2u \Delta u - 2\|\nabla u\|^2 = 2\langle \Delta v, v\rangle - 4u^2 -
2\|\nabla v\|^2$$ 
or
$$u \Delta u + 2u^2 = \langle \Delta v, v\rangle + \|\nabla u \|^2 - \|\nabla v\|^2.$$
We also know that $\|\nabla v\| \geq \|\nabla u\| = \|\nabla(\|v\|)\|$
and $-\langle \Delta v, v \rangle \geq -\|\Delta v\| \|v\| \geq
-bu$. Therefore 
$$-u(\Delta u + 2u) \geq -bu.$$
If $u \neq 0$ we have
$$-(\Delta u + 2u) \geq -b.$$
On the other hand if $u = 0$, $u$ has a local minimum so $-\Delta u
\geq 0$ and the inequality still holds. \qed{boundonLu} 

Define the operator $L$ by $L= -(\Delta + 2)$. Let $B_r$ be a ball of
radius $r$ centered at a point $p$. We
first need a fundamental solution for $L$. That is a radially
symmetric, smooth function $v(r)$ on $B_R \backslash p$ such that $Lv = 0$ and 
$$\int_{B_R} vL\phi dV = \lim_{\epsilon \rightarrow 0} \int_{B_R \backslash B_\epsilon}
vL\phi dV = \phi(p)$$ 
for all smooth functions $\phi$ with support in the interior of $B_R$.

\begin{lemma}
\label{fundamental}
The function
$$v(r) = \frac{-\cosh (\sqrt{3}r) + \coth (\sqrt{3}R) \sinh
  (\sqrt{3}r)}{4\pi \sinh r}$$ 
is a fundamental solution for $L$. Furthermore, $v(R) = 0$, $v(r) \leq
0 $ for $0 < r < R$ and 
$$v'(R) = \frac{\sqrt{3}}{4\pi \sinh(R) \sinh(\sqrt{3}R)}.$$
\end{lemma}

{\bf Proof.} For any radial function
$$-\Delta f  = \frac{\del^2 f}{dr^2} + 2\coth r \frac{df}{dr}.$$
Using this formula it is easy to check that $Lv = 0$.

The operator $L$ is self-adjoint so on any compact manifold $M$ with
boundary $\del M$, $L$ satisfies Green's identity: 
$$\int_M fLg dV = \int_M gLf dV + \int_{\del M} \left( f\frac{\del
    g}{\del n} - g\frac{\del f}{\del  n} \right)dA$$ 
where $f$ and $g$ are smooth functions on $M$ and $\frac{\del}{\del
  n}$ is the derivative in the direction of the outward normal. 

Applying Green's identity to $v$ and a test function $\phi$ we have:
\begin{eqnarray*}
\int_{B_R \backslash B_\epsilon} vL\phi dV & = & \int_{B_R \backslash B_\epsilon} \phi Lv dV +
\int_{\del B_R} \left( v\frac{\del \phi}{\del n} - \phi \frac{\del
    v}{\del r}\right) dA \\ & & - \int_{\del B_\epsilon} \left(
  v\frac{\del \phi}{\del n} - \phi \frac{\del v}{\del r}\right) dA\\ 
& = & \int_{\del B_\epsilon} v'\phi dA - \int_{\del B_\epsilon}
v\frac{\del \phi}{\del n} dA. 
\end{eqnarray*}

Clearly, $v(\epsilon) \sim -\frac{1}{4\pi \sinh R}$ for $\epsilon$ near $0$ and it
is easy to check that $v'(\epsilon) \sim \frac{1}{4\pi \sinh^2 r}$. We also
know that $\area (\del B_\epsilon) = 4\pi \sinh^2 (\epsilon)$ and
therefore 
$$\lim_{\epsilon \rightarrow 0} \int_{\del B_\epsilon} v' \phi dA =
\phi(p)$$ 
and
$$\lim_{\epsilon \rightarrow 0} \int_{\del B_\epsilon} v\frac{\del
  \phi}{\del n}dA = 0.$$ 
Therefore
$$\int_{B_R} vL\phi dA = \phi(p)$$
and $v$ is our fundamental solution.

The other properties of $v$ are a straightforward
calculation. \qed{fundamental}

\begin{lemma}
\label{meanbound}
Let $u$ be a smooth function on $B_R$ such that $Lu \geq -b$. Then
$$u(p) \leq \frac{1}{\sqrt{\vol(B_R)}} \sqrt{\int_{B_R} u^2 dV} + b/2.$$
\end{lemma}

{\bf Proof.} We first apply Green's identity on $B_{r'} \backslash B_\epsilon$:
\begin{eqnarray}
\int_{B_{r'}\backslash B_\epsilon} vLu dV & = & \int_{B_{r'}\backslash B_\epsilon} uLv dV + \int_{\del B_{r'}}
\left(v \frac{\del u}{\del r} - u \frac{\del v}{\del r} \right) dA
\nonumber \\ 
& & - \int_{\del B_\epsilon} \left(v\frac{\del u}{\del r} - u
  \frac{\del v}{\del r} \right) dA \label{green} 
\end{eqnarray}
where $v$ is the fundamental solution on the ball $B_{r'}$.
Recall that $Lv = 0$, $v(r') = 0$ and
$$\lim_{\epsilon \rightarrow 0} -\int_{\del B_\epsilon}
\left(v\frac{\del u}{\del r} - u \frac{\del v}{\del r} \right) dA =
u(p).$$ 

Therefore after taking the limit of (\ref{green}) as $\epsilon
\rightarrow 0$ and rearranging terms we have: 
$$u(p) = \int_{B_{r'}} v Lu dV + \int_{\del B_{r'}} u v' dA = \int_{B_{r'}} v Lu dV + v'(R)
\int_{\del B_{r'}}u dA.$$ 
By letting $u \equiv -1/2$ and solving the above equation for $\int_{B_{r'}}
vLu dV = \int_{B_{r'}} v dV$ we have 
$$0 \geq \int_{B_{r'}} v dV = \frac{\sqrt{3}}{2} \frac{\sinh r'}{\sinh (\sqrt{3}
  r')} - \frac{1}{2} \geq -\frac{1}{2}.$$ 
Now if $u$ is any smooth function on $B_{r'}$ with $Lu \geq -b$ we have:
$$\int_{B_{r'}} vLu dV \leq b/2.$$
Therefore:
$$u(p) \leq b/2 + \frac{\sqrt{3}}{4\pi\sinh(r') \sinh(\sqrt{3}r')}
\int_{\del B} udA.$$ 
Rearranging we have:
$$\frac{4\pi}{\sqrt{3}} (u(p) -b/2)(\sinh (\sqrt{3} r') \sinh r') \leq
\int_{\del B_{r'}} u dA.$$ 
Next we integrate both sides from $0$ to $R$:
$$\frac{4\pi}{\sqrt{3}} (u(p) -b/2) \int_0^R (\sinh (\sqrt{3} r') \sinh
r')dr \leq \int_0^R \left(\int_{\del B_{r'}} u dA\right) dr' = \int_{B_R} u
dV.$$ 
Since $\sinh r' \leq \sinh (\sqrt{3} r') / \sqrt{3}$ we have
$$\frac{4\pi}{\sqrt{3}} \int_0^R (\sinh (\sqrt{3} r') \sinh r') dr' \leq
4\pi \int_0^R \sinh^2 r' dr' = \vol (B_R).$$ 
By Holder's inequality
$$\int_{B_R} u dV \leq \sqrt{\vol(B_R)} \sqrt{\int_{B_R} u^2 dV}.$$
Therefore
$$u(p) \leq \frac{1}{\sqrt{\vol(B_R)}} \sqrt{\int_{B_R} u^2 dV} + b/2.$$
\qed{meanbound}

The three previous lemmas easily lead to the following theorem:
\begin{theorem}
\label{harmbound}
Let $v$ be a vector field on $B_R$ and assume $\|\Delta v\| < b$. Then: 
$$\|v(p)\| \leq \frac{1}{\sqrt{\vol (B_R)}}\sqrt{ \int_{B_R} \|v\|^2 dV} + b/2.$$
\end{theorem}

By Theorem 4.4 of \cite{Bromberg:rigidity} every cohomology class in $H^1(M; E)$ that
extends to a conformal deformation $\Phi$ of the projective boundary is represented by a Hodge form $\omega$. However, it is not shown that $\omega$ extends continuously to $\Phi$. We show this now.

\begin{theorem}
\label{context}
The Hodge form $\omega$ extends continuously to the holomorphic quadratic differential $\Phi$.
\end{theorem}

{\bf Proof.} We need to recall some of the proof of Theorem 4.4 in
\cite{Bromberg:rigidity}. In the proof $\omega = \omega_m + \omega_c$ where $\omega_m$ is a model deformation that extends continuously to $\Phi$ and the correction term $\omega_c$ is a trivial deformation. Then $\omega_m = ds_m$ where $s_m$ is the canonical lift of an automorphic vector field on $\tilde{M}$ while $\omega_c = ds_c$ where $s_c$ is the canonical lift of a vector field on $M$. The model $w_m$ is an
automorphic vector field that is in a standard form on the
geometrically finite ends and near the cone singularity. From this
standard form we know that the norms $\|\Delta
w_m\|$ and $\|\Delta (\curl w_m)\|$ decay to zero at the projective
boundary. 

Since we know that $\omega_m$ extends continuously to $\Phi$ to prove the theorem we need to show that $w_c$ and $\curl w_c$ extend continuously to the zero vector field. From the proof of Theorem 4.4 in \cite{Bromberg:rigidity} we know that
both $w_c$ and $\curl w_c$ are in $L^2$. We also know that $\Delta
w = \Delta (\curl w) = 0$ so $\Delta w_c = -\Delta w_m$ and $\Delta
(\curl w_c) = - \Delta (\curl w_m)$. Therefore the norms $\|\Delta
w_c\|$ and $\|\Delta (\curl w_c)\|$ decay to zero at the projective
boundary. 

We now apply Theorem
\ref{harmbound}. Let $p_n$ be a sequence of points in $M$ converging
to $p_\infty$ in $\Sigma$. For large values of $n$ there will be balls
$B_n$ centered at $p_n$ and embedded in $M$ such that $\|\Delta w_c\|$
and $\|\Delta (\curl w_c)\|$ is less than $b_n$ on $B_n$ with $b_n
\rightarrow 0$. Then by Theorem \ref{harmbound} 
$$\|w_c(p_n)\| \leq \frac{1}{\sqrt{\vol(B_n)}}\sqrt{\int_{B_n}
  \|w_c\|^2} + b_n/2$$ 
and
$$\|\curl w_c(p_n)\| \leq \frac{1}{\sqrt{\vol(B_n)}}\sqrt{ \int_{B_n}
  \|\curl w_c\|^2} + b_n/2.$$ 
Since both $w_c$ and $\curl w_c$ are in $L^2$ the right hand side of
these inequalities limits to zero at $n \rightarrow \infty$. Therefore
$\|w_c(p_n)\|$ and $\|\curl w_c(p_n)\|$ limit to zero on
$\Sigma$. \qed{context} 

\subsection{Strain fields}
Next we prove a mean value inequality for Hodge forms. The real and imaginary parts of a Hodge form are strain fields so this is equivalent to proving a mean value inequalities for harmonic strain fields.

We begin by defining a strain field. If $v$ is a vector field on a Riemannian manifold with covariant derivative $\nabla$ then $\nabla v$ is a vector value 1-form. The traceless symmetric part of $\nabla v$ is the {\em strain} $\str v$ of $v$ which measures the conformal distortion of $v$. The real and imaginary parts of an $E$-valued 1-form are vector valued 1-forms. It is shown in \cite{Hodgson:Kerckhoff:cone} that if $v$ is a harmonic, divergence free vector field and $\omega$ the associated Hodge form then $\Re \omega = \str v$ and $\Im \omega = -\str \curl v$.

\begin{lemma}
\label{harmonicstrain}
If $\eta$ is a stain field with $\Delta \eta = 0$ then $-\Delta \|\eta\| + 2\|\eta\| \geq 0$.
\end{lemma}

{\bf Proof.} The proof is a bit more complicated than Lemma \ref{boundonLu} because the Weitzenb\"{o}ck formula for strain fields is more involved. 
The Laplacian $\Delta = d\delta + \delta d$ for $E$-valued forms has the decomposition $\Delta = D^*D + DD^* + T^*T + TT^*$ where $d = D + T$ and $\delta = D^* + T^*$ are the decomposition of $d$ and $\delta$ into their real and imaginary parts. We deal with the first two terms $\Delta_D = D^*D + DD^*$ and last two terms $H = T^*T + TT^*$ separately. For $\Delta_D$ we have (see \cite{Wu:bochner})
$$\Delta_D \eta = -\underset{i}{\Sigma} \nabla^2_{e_i e_i} \eta - \underset{i,j}{\Sigma} \omega^i \wedge (R(e_i, e_j)\eta)(e_j).$$
Let
$$\cR \eta = - \underset{i,j}{\Sigma} \omega^i \wedge (R(e_i, e_j)\eta)(e_j)$$
so
$$\Delta_D = \nabla^* \nabla + \cR.$$
Since $R(e_i, e_j)$ is a tensor $\cR$ is purely algebraic and therefore easy to calculate.

Any strain field $\eta$ can be written as a linear combination $\eta = \underset{i,j}{\Sigma}f^l_k e_k \otimes \omega^l$. Then $\cR \eta = \underset{k, l}{\Sigma} f^l_k \cR(e_k \otimes \omega^l)$ so we need to calculate 
$$\cR(e_k \otimes \omega^l) = \underset{i,j}{\Sigma}\omega^i  \wedge \left(\omega^l(e_j) R(e_i, e_j)e_k \  + (R(e_i, e_j) \omega^l)(e_j) e_k\right).$$
To calculate these terms we recall that for hyperbolic space
$$R(e_i, e_j) e_k = \delta^i_k e_j - \delta^j_k e_i$$
where $\delta^i_j$ is the Kronecker delta function.
For the first term we have
$$\omega^i \wedge \omega^l(e_j) R(e_i, e_j)e_k = \delta^{l}_j(\delta^i_k e_j - \delta^j_k e_i) \otimes \omega^i.$$
If $k \neq l$ this is only non-zero if $i = j$ and $k = l$ in which case we get $e_l \otimes \omega^k$. If $k=l$ there are two non-zero terms $-e_i \otimes \omega^i$ when $i=j \neq k$. Therefore
$$\underset{i,j,k,l}{\Sigma} f^l_k \omega^i \wedge \omega^l(e_j) R(e_i, e_j)e_k =\eta^T - \tr \eta I.$$

For the second term we calculate
$$(R(e_i, e_j) \omega^l)(e_j) = \omega^l(-R(e_i, e_j) e_j) = (1-\delta^i_j) \omega^l(e_i) = (1-\delta^i_j) \delta^l_i$$
so
$$\underset{i,j}{\Sigma} \omega^i \wedge (R(e_i, e_j)\omega^l)(e_j) e_k = 2e_k \otimes \omega^l.$$

Since $\eta$ is a strain field it is traceless and symmetric therefore the two terms combine to give
$$\cR \eta = -3\eta.$$

For an harmonic strain field $\eta$ it is shown that $H\eta = \eta$ in \cite{Hodgson:Kerckhoff:cone}. Combining our work so far we have
$$\Delta \eta = \nabla^* \nabla \eta + \cR \eta + H \eta = \nabla^* \nabla \eta - 3\eta + \eta = 0$$
and therefore
$$\nabla^* \nabla \eta = 2\eta.$$

The remainder of the proof is exactly like the proof of Lemma \ref{boundonLu} and we will not repeat it. Note that for a harmonic vector field $\nabla^* \nabla v = -2v$ which accounts for the sign change from the bound we get for vector fields to the bound for strain fields.
\qed{harmonicstrain}

Now let $Lu = (-\Delta + 2)u$. We then can restate Lemmas \ref{fundamental} and \ref{meanbound} and Theorem \ref{harmbound} for this new definition of $L$. The proofs are so similar that we leave the details to the reader.

\begin{lemma}
\label{fundamentalstrain}
The function
$$v = \frac{-cos (\sqrt{2} r)  + \cot (\sqrt{2} R) \sin (\sqrt{2} r)}{4\pi \sinh r}$$
is a fundamental solution for $L$ if $R < \frac{\pi}{\sqrt{2}}$. Furthermore $v(R) = 0$, $v(r)\leq 0$ for $0< r<R$ and
$$v'(R) = \frac{\sqrt{2}}{4\pi \sinh (R) \sin (\sqrt{2} R)}.$$
\end{lemma}

\begin{lemma}
\label{meanboundstrain}
Let $u$ be a smooth function on $B_R$ such that $Lu \geq 0$. Then
$$u(p) \leq \frac{3 \sqrt{2\vol (B_R)}}{4\pi f(R)} \sqrt{\int_{B_R} u^2 dV}$$
where $f(R) = \cosh (R) \sin (\sqrt{2} R) - \sqrt{2} \sinh (R) \cos (\sqrt{2} R)$ for $R < \frac{\pi}{\sqrt{2}}$.
\end{lemma}

Using the fact that $\|\omega(p)\|^2 = \|\Re \omega (p)\|^2 + \|\Im \omega (p)\|^2$ we then have:

\begin{theorem}
\label{harmboundstrain}
Let $\omega$ be a Hodge form on a ball $B_R$ of radius $R$ centered at a point $p$. Then:
$$\|\omega(p)\| \leq \frac{3\sqrt{2 \vol(B_R)}}{4\pi f(R)} \sqrt{\int_{B_R} \|\omega\|^2 dV}$$
for $R < \frac{\pi}{\sqrt{2}}$.
\end{theorem}

\bibliographystyle{../tex/math}
\bibliography{../tex/math}

\begin{thebibliography}{CCHS}

\bibitem[And]{Anderson:projective}
C.~G. Anderson.
\newblock {Projective structures on Riemann surfaces and developing maps to
  $\hthree$ and $\cx P^n$}.
\newblock {\em Preprint} (1999).

\bibitem[AC1]{Anderson:Canary:cores}
J.~Anderson and R.~Canary.
\newblock {Cores of hyperbolic 3-manifolds and limits of Kleinian groups}.
\newblock {\em Amer. J. Math.} {\bf 118}(1996), 745--779.

\bibitem[AC2]{Anderson:Canary:corestwo}
J.~Anderson and R.~Canary.
\newblock {Cores of hyperbolic 3-manifolds and limits of Kleinian groups {\rm
  II}}.
\newblock {\em J. LMS} {\bf 61}(2000), 193--207.

\bibitem[BGS]{Ballman:Gromov:Schroeder}
W.~Ballman, M.~Gromov, and V.~Schroeder.
\newblock {\em Manifolds of Nonpositive Curvature}.
\newblock Birkhauser, 1985.

\bibitem[BB]{Brock:Bromberg:density}
J.~Brock and K.~Bromberg.
\newblock {On the density of geometrically finite Kleinian groups}.
\newblock {\em In preparation} (2002).

\bibitem[BBES]{Brock:Bromberg:Evans:Souto}
J.~Brock, K.~Bromberg, R.~Evans, and J.~Souto.
\newblock {Boundaries of deformation spaces and Ahlfors' measure conjecture}.
\newblock {\em 2002 Preprint available at
  \verb+front.math.ucdavis.edu/math.GT/0211022+}.

\bibitem[Br1]{Bromberg:projective}
K.~Bromberg.
\newblock {Projective structures with degenerate holonomy and the Bers' density
  conjecture}.
\newblock {\em 2002 Preprint available at
  \verb+www.math.caltech.edu/people/bromberg.html+}.

\bibitem[Br2]{Bromberg:rigidity}
K.~Bromberg.
\newblock {Rigidity of geometrically finite hyperbolic cone-manifolds}.
\newblock {\em 2002 Preprint available at
  \verb+www.math.caltech.edu/people/bromberg.html+}.

\bibitem[Can]{Canary:retraction}
R.~D. Canary.
\newblock {The conformal boundary and the boundary of the convex core}.
\newblock {\em Duke Math. J.} {\bf 106}(2000), 193--207.

\bibitem[CCHS]{Canary:Culler:Hersonsky:Shalen}
R.~D. Canary, M.~Culler, S.~Hersonsky, and P.~B. Shalen.
\newblock {Approximation by maximal cusps in the boundaries of quasiconformal
  deformation spaces}.
\newblock {\em Preprint} (2000).

\bibitem[CEG]{Canary:Epstein:Green}
R.~D. Canary, D.~B.~A. Epstein, and P.~Green.
\newblock {Notes on notes of Thurston}.
\newblock In {\em Analytical and Geometric Aspects of Hyperbolic Space}, pages
  3--92. Cambridge University Press, 1987.

\bibitem[CM]{Canary:Minsky:tameness}
R.~D. Canary and Y.~N. Minsky.
\newblock {On limits of tame hyperbolic 3-manifolds}.
\newblock {\em J. Diff. Geom.} {\bf 43}(1996), 1--41.

\bibitem[Ep]{Epstein:horospheres}
C.~Epstein.
\newblock {Envelopes of horospheres and Weingarten surfaces in hyperbolic
  3-space}.
\newblock {\em preprint}.

\bibitem[Ev]{Evans:tameness}
R.~Evans.
\newblock {Tameness persists}.
\newblock {\em preprint} (2001).

\bibitem[HK1]{Hodgson:Kerckhoff:cone}
C.~Hodgson and S.~Kerckhoff.
\newblock {Rigidity of hyperbolic cone-manifolds and hyperbolic Dehn surgery}.
\newblock {\em J. Diff. Geom.} {\bf 48}(1998), 1--59.

\bibitem[HK2]{Hodgson:Kerckhoff:harmonic}
C.~Hodgson and S.~Kerckhoff.
\newblock {Harmonic deformations of hyperbolic 3-manifolds}.
\newblock {\em 2002 Preprint}.

\bibitem[HK3]{Hodgson:Kerckhoff:tube}
C.~Hodgson and S.~Kerckhoff.
\newblock {The shape of hyperbolic Dehn surgery space}.
\newblock {\em In preparation}.

\bibitem[HK4]{Hodgson:Kerckhoff:dehn}
C.~Hodgson and S.~Kerckhoff.
\newblock {Universal bounds for hyperbolic Dehn surgery}.
\newblock {\em 2002 Preprint available at
  \verb+front.math.ucdavis.edu/math.GT/0204345+}.

\bibitem[Ko1]{Kodani:convergence}
S.~Kodani.
\newblock {Convergence theorem for Riemannian manifolds with boundary}.
\newblock {\em Compositio Math.} {\bf 75}(1990), 171--192.

\bibitem[Ko2]{Kojima:cone}
S.~Kojima.
\newblock {Deformations of hyperbolic 3-cone manifolds}.
\newblock {\em J. Diff. Geom.} {\bf 49}(1998), 469--516.

\bibitem[Mc]{McMullen:cusps}
C.~McMullen.
\newblock {Cusps are dense}.
\newblock {\em Annals of Math.} {\bf 133}(1991), 217--247.

\bibitem[Ot]{Otal:unknottwo}
J.~P. Otal.
\newblock {{Les g\'eod\'esiques ferm\'ees d'une vari\'et\'e hyperbolique en
  tant que noeuds}}.
\newblock {\em Preprint} (2002).

\bibitem[RS]{Rivin:Schlenker:schlafli}
I.~Rivin and J-M. Schlenker.
\newblock {{On the Schl\"afli differential formula}}.
\newblock {\em preprint} (1998).

\bibitem[Wu]{Wu:bochner}
H.~Wu.
\newblock {The Bochner technique in differential geometry}.
\newblock {\em Math. Reports, London} (1987).

\end{thebibliography}

\begin{sc}
\noindent
Department of Mathematics\\
California Institute of Technology\\
Mailcode 253-37\\
Pasadena CA 91125\\
\end{sc}

\end{document}